\documentclass[sn-mathphys,Numbered]{sn-jnl}

\usepackage{amsfonts,amsmath,amssymb,amsthm}
\usepackage[export]{adjustbox}
\usepackage{mathrsfs}
\usepackage{fixmath}
\usepackage{graphicx}%
\usepackage{multirow}%
\usepackage{caption}
\usepackage{subcaption}
\usepackage{float}
\usepackage{tikz}
\usepackage[figuresright]{rotating}
\usepackage{array}
\usepackage{amsthm}%
\usepackage{mathrsfs}%
\usepackage[title]{appendix}%
\usepackage{xcolor}%
\usepackage{textcomp}%
\usepackage{manyfoot}%
\usepackage{booktabs}%
\usepackage{algorithm}%
\usepackage{algorithmicx}%
\usepackage{algpseudocode}%
\usepackage{tabularx}
\usepackage{listings}%

\newcommand{\bpsi}{\boldsymbol{\psi}}

\newcommand{\ds}{\displaystyle}

\begin{document}

\title[An Exponentiated Odds Ratio Generator]{\vspace{-2cm} 
\textbf{Advancing Continuous Distribution Generation: An Exponentiated Odds Ratio Generator Approach}}


\author[1]{\fnm{Xinyu} \sur{Chen}}
\equalcont{These authors contributed equally to this work.}

\author[2]{\fnm{Yuanqi} \sur{Xie}}

\author[1]{\fnm{Achraf} \sur{Cohen}}

\author*[1]{\fnm{Shusen} \sur{Pu}}\email{spu@uwf.edu}
\equalcont{These authors contributed equally to this work.}

\affil*[1]{\orgdiv{Department of Mathematics and Statistics}, \orgname{University of West Florida}, \orgaddress{\street{11000 University Pkwy}, \city{Pensacola}, \postcode{32514}, \state{FL}, \country{USA}}}

\affil[2]{\orgdiv{Department of Computer Science}, \orgname{University of West Florida}, \orgaddress{\street{11000 University Pkwy}, \city{Pensacola}, \postcode{32514}, \state{FL}, \country{USA}}}



\abstract{This paper presents a new methodology for generating continuous statistical distributions, integrating the exponentiated odds ratio within the framework of survival analysis. This new method enhances the flexibility and adaptability of distribution models to effectively address the complexities inherent in contemporary datasets. The core of this advancement is illustrated by introducing a particular subfamily, the ``Type-2 Gumbel Weibull-G Family of Distributions." We provide a comprehensive analysis of the mathematical properties of these distributions, encompassing statistical properties such as density functions, moments, hazard rate and quantile functions, R\'enyi entropy, order statistics, and the concept of stochastic ordering. To establish the robustness of our approach, we apply five distinct methods for parameter estimation. The practical applicability of the Type-2 Gumbel Weibull-G distributions is further supported through the analysis of three real-world datasets. These empirical applications illustrate the exceptional statistical precision of our distributions compared to existing models, thereby reinforcing their significant value in both theoretical and practical statistical applications.}




\maketitle

\section{Introduction}

With the ever-increasing complexity and volume of data across various disciplines, developing and applying a diverse array of statistical distributions has become a paramount area of research. These continuous distributions are essential in modeling, forecasting, and interpreting complex data, facilitating hidden patterns and relationships \cite{land12010035}. Nonetheless, the rapid growth and evolving nature of contemporary data often pose challenges that traditional distributions struggle to cope with, leading to a need for novel statistical distribution generators \cite{JIANG2023102309}.

Over the past several decades, numerous methods for creating new continuous probability distributions have been explored. These include methods such as transformations of random variables, the use of mixed models, and advanced compounding methods \cite{cooray2006generalization,bourguignon2014weibull, Pu2016GEMW, oluyede2018gamma}. Among the recent contributions in this field are the gamma-Topp-Leone-Type II-exponentiated half Logistic-G distribution family \cite{oluyede2023gamma}, an enhanced version of the generalized Weibull distribution \cite{shama2023modified}, the innovative modified alpha power Weibull-X distribution set \cite{emam2023modeling}, the Topp-Leone type II exponentiated half logistic-G distribution category \cite{gabanakgosi2023topp}, the modified-half-Normal distribution by \cite{sun2023modified}, an inventive extension of the power Lindley distribution \cite{sun2023modified}, the truncated inverse generalized Rayleigh distribution \cite{guptha2023new}, the Ristić–Balakrishnan–Topp–Leone–Gompertz-G family of distributions \cite{Pu2023JSTA}, the odd Gompertz-G distribution family \cite{kajuru2023odd}, the inverse Burr-generalized distribution series \cite{osagie2023inverse}, the gamma inverse paralogistic distribution \cite{marasigan2023new}, and the shifted generalized truncated Nadarajah-Haghighi distribution by \cite{azimi2023new}.

This paper presents a new generator of continuous statistical distributions, designed to offer greater flexibility and adaptability in addressing the increasing complexity of contemporary datasets. The foundation of this generator's development lies in survival analysis, drawing inspiration from the extreme values of odds ratio data, as elaborated in Section \ref{sec: Newfamily}. The generator proposed herein aims to provide a comprehensive framework for the creation of a broad range of distributions, each characterized by unique shapes and properties. This feature allows for more precise modeling and analysis of an extensive variety of data. The adaptability of the proposed generator is illustrated through its application in deriving new families of distributions, demonstrating superior data fitting capabilities compared with traditional methods.

The structure of the remaining sections of this paper is outlined as follows: Section 2 provides a description of the methodology of the new generator. Section 3 presents the new subfamily of distributions generated by the proposed technique and investigates the mathematical properties of the resulting distributions. Four estimation methods are discussed in Section 4 to test the robustness of the model.
Section 5 introduces several special cases of the new family of distributions. In Section 6, we evaluate the performance and flexibility of the generator through a set of experiments and comparative analyses. Finally, we conclude the paper and discuss the future directions in Section 7.

\section{The New Generator Based on Odds Ratio}\label{sec:distribution}

The field of survival analysis is a critical component in various research domains, such as medical science, engineering, and social sciences. Central to survival analysis is the odds ratio, which quantifies the probability that an individual or component, defined by a specific lifespan following a continuous distribution $H(x,\bpsi)$, will fail or expire at a particular point in time, $x$ \cite{kleinbaum1996survival}. The odds ratio, expressed as $\ds {H(x,\bpsi)}/{\overline{H}(x,\bpsi)}$, has become an important tool for understanding and assessing risk factors, shedding light on the relative probabilities of event outcomes, namely, death or failure \cite{tsiatis2023estimation, vanderweele2020optimal, penner2019odds}.

This study presents a novel generator developed by integrating the methodologies delineated in \cite{bourguignon2014weibull} and \cite{Alzaatreh2013}. This new generator can be efficiently applied to any baseline distribution with a cumulative distribution function (cdf), represented as $R(t)$. It is mathematically defined by:
\begin{align}
F_{{\text{RT--EOR-H}}}(x;\alpha,\beta)&= \int_{0}^{-\log (1-F_{WG}(x))}r(t) dt \\
&=R \left\{-\log (1-F_{WG}(x)) \right\} \\
&=R\left\{\alpha \left[ \frac{H(x,\bpsi)}{\overline{H}(x,\bpsi)} \right]^\beta \right\} \label{eq:new_Gx}
\end{align}
Here, $r(t)$ and $R(t)$ denote the probability density function (pdf) and cdf of a continuous random variable, respectively. $f_{WG}(x)$, and $F_{WG}(x)$ represent the pdf and cdf of the Weibull-G family of distributions proposed by Bourguignon et al. \cite{bourguignon2014weibull}. $H(x,\bpsi)$ and $\overline{H}(x,\bpsi)$ represent the cdf and the survival function of the baseline distribution, respectively. This newly designed generator of continuous distributions is denoted as the R-transferred exponentiated odds ratio generator (RT-EOR-H), denoted as $F_{_{\text{RT--EOR-H}}}(x;\alpha,\beta)$.

This new approach provides a comprehensive framework, allowing researchers to delve into the exponentiated odds ratio of any baseline function $H(x,\bpsi)$, embedded within any conceivable distribution $R(t)$.  Crucially, this new generator, by permitting the investigation of any baseline function within any distribution $R$, significantly enhances the flexibility of survival analyses, yielding more accurate, adaptable, and precise predictions. These advancements meet the need for continued evolution in survival analysis, providing a path for enhanced comprehension of the complexities of survival and risk.

Some example subfamilies are listed in Table \ref{functionslist}. 

\begin{table}[htp]
    \centering
    \caption{Families of distributions derived from different $R(t)$}
    \label{functionslist}
    \setlength{\tabcolsep}{2mm}
    \begin{tabular}{lll}
    \hline
	Distribution &  $R(t)$ & RT-EOR-H \\
	\hline
        Uniform & $\frac{t}{\theta}$ & $\frac{\alpha}{\theta} \left[\frac{H(x,\bpsi)}{\overline{H}(x,\bpsi)}\right]^\beta$ \\
        Normal&$\Phi\left(\frac{x-\mu}{\sigma}\right)$ & $\Phi \left(\frac{\alpha \left[ \frac{H(x,\bpsi)}{\overline{H}(x,\bpsi)} \right]^\beta-\mu}{\sigma} \right)$ \\
        Gamma & $\frac{1}{\Gamma(k)} \gamma (k,\frac{t}{\theta})$ & $\frac{1}{\Gamma(k)}\gamma \left(k,\frac{\alpha}{\theta}\left[\frac{H(x,\bpsi)}{\overline{H}(x,\bpsi)}\right]^\beta \right)$\\
        Log-logistic& $\frac{1}{1+(t/c)^{-k}}$ & $\frac{1}{1+\left(\frac{\alpha}{c} \left[ \frac{H(x,\bpsi)}{\overline{H}(x,\bpsi)} \right]^\beta\right)^{-k}}$\\
        Rayleigh & $1-e^{-\frac{t^2}{2\sigma^2}}$ &$1-\exp\left\{-\frac{\alpha^2}{2\sigma^2} \left[ \frac{H(x,\bpsi)}{\overline{H}(x,\bpsi)} \right]^{2\beta}\right\}$\\
        Weibull & $1-e^{-(t/\lambda)^k}$ & $1-\exp\left\{-\left(\frac{\alpha}{\lambda} \left[ \frac{H(x,\bpsi)}{\overline{H}(x,\bpsi)} \right]^\beta\right)^k \right\}$\\
        Type-2 Gumbel & $e^{-\lambda t^{-\delta}}$ & $\exp\left\{-\lambda  \left(\alpha\left[ \frac{H(x,\bpsi)}{\overline{H}(x,\bpsi)} \right]^\beta\right)^{-\delta}\right\}$\\
        Lomax&$1-(1+\lambda t)^{-k}$ & $1-\left(1+\lambda \alpha\left[ \frac{H(x,\bpsi)}{\overline{H}(x,\bpsi)} \right]^\beta \right)^{-k}$\\
        Burr XII & $1-(1+t^c)^{-k}$ & $1-\left\{1+\left(\alpha\left[ \frac{H(x,\bpsi)}{\overline{H}(x,\bpsi)} \right]^\beta\right)^c \right\}^{-k}$\\
        Pareto & $1-\left(\frac{\gamma}{t}\right)^k$ & $1-\left(\frac{\gamma}{\alpha\left[ \frac{H(x,\bpsi)}{\overline{H}(x,\bpsi)} \right]^\beta}\right)^k $\\
        L\'evy & $\text{erfc}\left( \sqrt{\frac{c}{2t}}\right)$ & $\text{erfc}\left( \sqrt{\frac{c}{2\alpha\left[ \frac{H(x,\bpsi)}{\overline{H}(x,\bpsi)} \right]^\beta}}\right)$ \\
        Fr\'echet & $e^{-\lambda t^\gamma}$ & $\exp \left\{-\lambda \left(\alpha\left[ \frac{H(x,\bpsi)}{\overline{H}(x,\bpsi)} \right]^\beta\right) ^\gamma \right\}$\\
        Kumaraswamy & $(1-x^k)^{-\lambda}-1$ & $\left\{ 1-\left(\alpha\left[ \frac{H(x,\bpsi)}{\overline{H}(x,\bpsi)} \right]^\beta\right)^k\right\}^{-\lambda} - 1$\\
	\hline
	\end{tabular}
\end{table}

\section{The Type-2 Gumbel Weibull-G Family of Distributions}
\label{sec: Newfamily}
To illustrate the applicability of the new generator, we will focus on one of its subfamilies in this paper. Here, we combine the generator with the cdf of Type-2 Gumbel distribution as follows
\begin{align}
    R(t)= e^{-\lambda t^{-\delta}} 
    \label{cdf_R1}
\end{align}
Then we obtain the cdf and pdf of a new Type-2 Gumbel Weibll-G (T2GWG)family of distributions as
\begin{align}
    F_{_{T2GWG}}(x) = \exp\left\{-\lambda\alpha^{-\delta} \left[ \frac{H(x,\bpsi)}{\overline{H}(x,\bpsi)} \right]^{-\beta \delta} \right\}
    \label{cdf}
\end{align}
To avoid over-parameterization, we let $\lambda=1$,  $\delta=1$, and substitite $\alpha=\alpha^{-1}$, now the cdf reduces to
\begin{align}
    F_{_{T2GWG}}(x) &= \exp\left\{-\alpha \left[ \frac{H(x,\bpsi)}{\overline{H}(x,\bpsi)} \right]^{-\beta } \right\}
    \label{cdf_reduce}
\end{align}
and 
\begin{align}
    f_{_{T2GWG}}(x) = \alpha \beta h(x,\bpsi) \frac{H(x,\bpsi)^{-\beta -1}}{\overline{H}(x,\bpsi)^{-\beta  +1 }} \exp\left\{-\alpha \left[ \frac{H(x,\bpsi)}{\overline{H}(x,\bpsi)} \right]^{-\beta }\right\}
    \label{pdf}
\end{align}


An interpretation of the T2GWG family of distributions can be given as follows.
Let $Z$ be a lifetime random variable with a baseline cdf $H(x,\bpsi)$. The odds ratio that an individual following the lifetime $Z$ will die (failure) at time $x$
is $\ds \frac{H(x,\bpsi)}{\overline{H}(x,\bpsi)}$. For a sequence of such independent and identical odds ratios $Z_1, Z_2, \cdots, Z_n$, the maximum $M_n=\text{max}(Z_1, Z_2, \cdots, Z_n)$, then $M_n$ follows a distribution which converges to the T2GWG distribution as provided in equation\eqref{cdf_reduce} \cite{cooray2006generalization,bourguignon2014weibull}.
 
In the subsequent subsections, we delve into the statistical properties of this novel distribution. Our discussion will cover a broad range of topics, including the expansion of the density function, hazard rate and quantile functions, moments, incomplete moments, the generation function, R\'enyi Entropy, order statistics, and stochastic ordering.

\subsection{Expansion of the pdf}
Consider the following expansion
\begin{align}
    \exp\left\{-\alpha \left[ \frac{H(x,\bpsi)}{\overline{H}(x,\bpsi)} \right]^{-\beta }\right\} = \sum_{j=0}^{\infty} \frac{(-1)^j \alpha^j}{j!} \left[ \frac{H(x,\bpsi)}{\overline{H}(x,\bpsi)} \right]^{-j\beta},
\end{align}
then the pdf of T2GWG can be expanded as 
\begin{align}
    f_{_{T2GWG}}(x)& = \alpha \beta h(x,\bpsi) \frac{H(x,\bpsi)^{-\beta-1}}{\overline{H}(x,\bpsi)^{-\beta+1}}\sum_{j=0}^{\infty}\frac{(-1)^j\alpha^j}{j!}\left[ \frac{H(x,\bpsi)}{\overline{H}(x,\bpsi)} \right]^{-j\beta} \nonumber \\
    & = \alpha \beta h(x,\bpsi) \sum_{j=0}^{\infty}\frac{(-1)^j\alpha^j}{j!}\frac{\overline{H}(x,\bpsi)^{\beta+j\beta-1}}{H(x,\bpsi)^{\beta+j\beta+1}}
\end{align}
Moreover, note that 
\begin{align}
    \overline{H}(x,\bpsi)^{\beta +j\beta-1} & =[1-H(x,\bpsi)]^{\beta+j\beta-1} \nonumber \\
    &=\sum_{k=0}^{\infty} \binom{\beta(j+1)-1}{k}
    (-1)^k H(x,\bpsi)^k
\end{align}
Then
\begin{align}
    f_{_{T2GWG}}(x)& = \alpha \beta h(x,\bpsi)\sum_{j=0}^{\infty}\frac{(-1)^j \alpha^j}{j!} \sum_{k=0}^{\infty}
    \binom{\beta(j+1)-1}{k}
    \frac{(-1)^k H(x,\bpsi)^k}{ H(x,\bpsi)^{\beta+j\beta+1}} \nonumber \\
    &=\alpha \beta h(x,\bpsi) \sum_{j,k=0}^{\infty}\frac{(-1)^{j+k} \alpha^j}{j!}\binom{\beta(j+1) -1}{k}\left[H(x,\bpsi) \right]^{k-\beta-j\beta-1}\nonumber \\
    &=\sum_{j,k=0}^{\infty} c_{j,k} r_{k-\beta(j+1) -1}(x,\bpsi)
\end{align}
where 
\begin{align}
    c_{j,k}=\frac{\alpha\beta(-1)^{j+k}
    \alpha^j}{j!(k-\beta(j+1))} \binom{\beta (j+1)-1}{k} \label{cij}
\end{align}
and
$$r_{k-\beta-j\beta-1}(x,\bpsi)=(k-\beta(j+1)) h(x,\bpsi)H(x,\bpsi)^{k-\beta(j+1)-1}$$
which is the pdf of exponentiated generalized (EG) distribution with parameter $\beta^*=k-\beta (j+1)$.

\subsection{Hazard rate and quantile functions}
Building upon our earlier discussion on the odds ratio and survival analysis, we now delve into the key mathematical structures: the hazard rate function and the quantile function. These are crucial to the deeper understanding and application of survival analysis as they allow the computation of survival probabilities and survival times.

\subsubsection{Hazard rate and quantile functions}
The hazard rate function (hrf) plays an important role in survival analysis as it defines the instantaneous potential per unit time for the occurrence of an event given survival up to that time. On the other hand, the quantile function is essential in determining the time at which a certain proportion of survival is expected. In the newly proposed family of distributions, these functions take a particularly interesting form as follows.
\begin{align}
    h_{_{T2GWG}}(x)=\frac{\alpha \beta h(x,\bpsi) \frac{H(x,\bpsi)^{-\beta-1}}{\overline{H}(x,\bpsi)^{-\beta+1}} \exp\left\{-\alpha\left[ \frac{H(x,\bpsi)}{\overline{H}(x,\bpsi)} \right]^{-\beta} \right\}}{1-\exp\left\{-\alpha \left[\frac{H(x,\bpsi)}{\overline{H}(x,\bpsi)} \right]^{-\beta} \right\}}
\end{align}
In addition to the hrf, it is also useful to consider its reciprocal, termed the reverse hazard rate function. This function essentially reflects the hazard function's properties but is viewed from the perspective of the event not occurring.

\begin{align}
    \tau_{_{T2GWG}}(x)&= \frac{\alpha \beta h(x,\bpsi) \frac{H(x, \bpsi)^{-\beta -1}}{\overline{H}(x, \bpsi)^{-\beta +1}}\exp\left\{-\alpha \left[ \frac{H(x,\bpsi)}{\overline{H}(x, \bpsi)} \right]^{-\beta} \right\}}{\exp\left\{-\alpha \left[\frac{H(x, \bpsi)}{\overline{H}(x,\bpsi)} \right]^{-\beta } \right\}} \nonumber \\
    &=\alpha \beta h(x,\bpsi) \frac{H(x,\bpsi)^{-\beta-1}}{\overline{H}(x,\bpsi)^{-\beta+1}}
\end{align}
\subsubsection{Quantile function}
Next, we will focus on the quantile function, which is particularly useful when determining the survival time corresponding to a specific survival probability.
$$F_{_{T2GWG}}(x)= \exp\left\{-\alpha \left[ \frac{H(x,\bpsi)}{\overline{H}(x,\bpsi)}\right]^{-\beta}\right\}=p$$
for $0\leq p\leq 1$, 
then, it is sufficient to solve 
\begin{align}
    H(x,\bpsi)= \displaystyle \frac{1}{\left(\frac{\log p}{-\alpha} \right)^{\frac{1}{\beta}}+1} :=q
\end{align}
Thus, the quantile $x_p$ of the distribution reduces to the quantile $x_q$ of the baseline distribution with cdf $H(x,\bpsi)$ and is given by 
\begin{align}
    x_q = H^{-1}(q)
\end{align}

\subsection{Moments, incomplete moments and generating functions}
\subsubsection{Moments}
In the field of statistics, moments play a vital role in characterizing the properties of a probability distribution. Moments provide important summary measures of the characteristics of data sets. The first moment about the origin, also known as the mean, measures the location of the distribution. The second central moment is known as the variance, which quantifies the spread or dispersion of the distribution. The third and fourth moments, skewness and kurtosis, respectively, describe the shape of the distribution, capturing aspects of its asymmetry and tailedness.
We can present the $r^{th}$ moment of the distributions as
\begin{align}
    E(Y^r)=\sum_{j,k=0}^{\infty}c_{j,k}
    E(Z_{j,k}^{r})
\end{align}
where $Z_{j,k}$ is the exponentiated generalized distribution with the parameter $\beta^*=k-\beta (j+1)$ and $c_{j,k}$ is defined by Eq. \ref{cij}. 

\subsubsection{Incomplete Moments, Conditional Moments and Moment Generating Function}
While moments give us an understanding of the general characteristics of a distribution, they do not always provide sufficient detail about specific intervals or subsets within the data. This is where incomplete moments come into play. Incomplete moments, also known as truncated or restricted moments, are defined similarly to regular moments but are integrated over a subset of the possible range of the variable. They offer a more granular insight into the characteristics of the distribution within specific ranges. This makes them particularly useful when analyzing left or right-skewed data or when assessing the impact of outlier observations.
The incomplete moment is provided as
\begin{align}
    I_{Y}(z)=\int_0^z y^s f_{_{T2GWG}}(y) dy = \sum_{j,k=0}^{\infty} c_{j,k} I_{j,k}(y)
\end{align}
where $I_{j,k}(y)=\int_0^z y^s r_{k-\beta (j+1)}(x,\bpsi)$. \\
The $rth$ conditional moments of the Type-2 Gumble Weibull-G family of distributions is given by 
\begin{align}
    E(Y^r|Y\geq a) & = \frac{1}{\overline{F}_{T2GWG}(a;\alpha,\beta,\bpsi)} \int_{t}^{\infty} y^r f_{_{T2GWG}} (y;\alpha,\beta,\bpsi) dy \nonumber \\
    & = \frac{1}{\overline{F}_{T2GWG}(a;\alpha,\beta,\bpsi)} \sum_{j,k=0}^\infty c_{j,k} I_{j,k}(y)
\end{align}
where $I_{j,k}(y)$ is defined above. \\
The moment generating function is given by 
\begin{align}
    M_{Y}(t) = E(e^{tY}) = \sum_{j,k=0}^\infty c_{j,k} E(e^{tZ_{j,k}}) = \sum_{j,k=0}^\infty c_{j,k} M_{Z_{j,k}}(t)
\end{align}
and the characteristic function is defined as 
\begin{align}
    \phi(t) = E(e^{itY}) = \sum_{j,k=0}^\infty c_{j,k} E(e^{itZ_{j,k}}) = \sum_{j,k=0}^\infty c_{j,k} \phi_{k-\beta (j+1)}(t),
\end{align}
where $\phi_{k-\beta (j+1)}(t)$ is the characteristic function of EG distribution with parameter $\beta^*=k-\beta (j+1)$.

\subsection{R\'enyi Entropy}
Rényi entropy is named after the Hungarian mathematician Alfr\'ed R\'enyi \cite{Renyi1961Measures}. While Shannon entropy is perhaps the most commonly referenced form of entropy in the field of information theory, characterizing the average uncertainty or unpredictability of a source of information, Rényi entropy provides a more generalized measure. Rényi entropy finds applications in various domains including physics, computer science, statistics, and quantum information theory. For example, in the context of machine learning, it can be used to measure the diversity or complexity of learned models. In statistical physics, it is useful in understanding the thermodynamics of complex systems.
The R\'enyi Entropy of this new distribution is 
\begin{align}
    I_R(\omega)=&(1-\omega)^{-1} \log\left[\int_{-\infty}^{\infty} f^{\omega}(x)dx\right] \nonumber \\
    =&(1-\omega)^{-1}\left\{\omega( \log\alpha+\log\beta) \right.\nonumber\\ 
    &+\left.\log\left[\int_{-\infty}^{\infty} h^\omega(x,\bpsi) \frac{H(x,\bpsi)^{\omega(-\beta -1)}}{\overline{H}(x,\bpsi)^{\omega(-\beta+1)}}\exp \left\{-\omega \alpha \left[ \frac{H(x,\bpsi)}{\overline{H}(x,\bpsi)} \right]^{-\beta} \right\} \right] \right\}  \nonumber
\end{align}
where $\omega>0$ and $\omega \neq 1$. 
Apply the same expansion technique for the pdf, we obtain
\begin{align}
    I_R(\omega)&=(1-\omega)^{-1} \left\{\omega(\log\alpha + \log\beta) \right. \nonumber \\ 
    &+\left.\log \left[ \sum_{i=0}^{\infty} \frac{(-1)^i(\omega \alpha)^i}{i!}\int_{-\infty}^{ \infty} h^{\omega}(x,\bpsi)\frac{ [\overline{H}(x,\bpsi)]^{\omega( \beta-1)+ i\beta}}{[H(x,\bpsi)]^{\omega(\beta+1)+ i\beta}} dx\right]\right\}
\end{align}
Consider that
\begin{align}
    \overline{H}(x,\bpsi)^{\omega(\beta-1)+ i\beta}
    & = [1-H(x,\bpsi)]^{\omega(\beta -1)+i\beta}\nonumber \\ 
    & = \sum_{j=0}^{\infty} \binom{\omega(\beta-1)+ i\beta}{j} (-1)^j H(x,\bpsi)^j
\end{align}
Thus we can write the Re\'nyi Entropy as 
\begin{align}
    I_R(\omega) =&(1-\omega)^{-1} \left\{\omega(\log\alpha+\log\beta)
    + \log\left[\sum_{i=0}^{\infty} \sum_{j=0}^{\infty}\frac{(-1)^{i+j}(\omega \alpha)^i}{i!}\right.\right. \nonumber \\ 
    & \times \left. \left. \binom{\omega(\beta-1)+i\beta}{j} 
    \int_{-\infty}^{\infty}h^{\omega}(x, \bpsi)(H(x,\bpsi))^{j-\omega (\beta+1)-i\beta}dx\right]\right\} \nonumber \\
    = & (1-\omega)^{-1} \left\{\omega (\log\alpha+\log\beta) +\log\left[ \sum_{i=0}^{\infty}\sum_{j=0}^{\infty}\frac{(-1)^{i+j}(\omega \alpha)^i}{i!} \right. \right. \nonumber \\
    & \times \left. \binom{\omega (\beta-1)+i\beta}{j} \frac{\omega^\omega}{[j-\omega(\beta+1)-i\beta+\omega ]^\omega}\right.  \nonumber \\ 
    & \times \left. \left. \int_{-\infty}^{\infty} [\frac{j-\omega (\beta+1)-i\beta+\omega}{\omega} h(x,\bpsi)(H(x,\bpsi))^{\frac{j-\omega(\beta+1)-i\beta}{\omega}}]^{\omega} \right] \right\} \nonumber \\
    = & (1-\omega)^{-1} \left\{\omega (\log\alpha+\log\beta) + \log \left[\sum_{i=0}^{\infty} \sum_{j=0}^{\infty} 
    \frac{(-1)^{i+j}(\omega \alpha)^i}{i!}  \right. \right. \nonumber \\
    & \times \binom{\omega(\beta-1)+i\beta}{j}\frac{\omega^\omega}{[j-\omega(\beta+1)-i\beta+ \omega]^\omega}\times e^{(1-\omega)I_{REG}} 
\end{align}
where $I_{REG}$ is the R\'enyi entropy of the exponentiated generalized distribution with parameter $\ds \beta^*=\frac{j-\omega(\beta+1)-i\beta+\omega}{\omega}$.

\subsection{Order statistics}
In the field of statistics, order statistics are a fundamental concept that allows for deeper analysis and understanding of sampled data. Specifically, order statistics are the values from a random sample sorted in ascending or descending order. This sorting process provides a powerful perspective on the sample's overall distribution and associated characteristics.
Order statistics are used in a variety of applications, including non-parametric statistics (which does not rely on parameters defined in terms of a theoretical or assumed population), reliability engineering, and statistical quality control. They also play a central role in the construction of quantile-quantile plots, which are used to assess if a data set follows a particular theoretical distribution.

Let $X_1, X_2,...,X_n$ be independent identically distributed random variables distributed by Eqn.~\eqref{pdf}. The pdf of the $i^{th}$ order statistic $f_{i:n}(x)$ is given by
\begin{align}
   f_{i:n}(x) & = \frac{n!f_{_{T2GWG}}(x)}{(i-1)!(n-i)!}[F_{_{T2GWG}}(x)]^{i-1}[1-F_{_{T2GWG}}(x)]^{n-i} \nonumber \\
   & = \frac{n!f_{_{T2GWG}}(x)}{(i-1)!(n-i)!} \sum_{m=0}^{n-i} \binom{n-i}{m} (-1)^m[F_{_{T2GWG}}(x)]^{i-1+m}  \nonumber \\
   & = \frac{n!f_{_{T2GWG}}(x)}{(i-1)!(n-i)!} \sum_{m=0}^{n-i} \binom{n-i}{m}(-1)^m \exp\left\{(i-1+m)(-\alpha) \left[\frac{H(x,\bpsi)}{\overline{H}(x,\bpsi)} \right] ^{-\beta}
   \right\} \nonumber \\
   & = \frac{n!}{(i-1)!(n-i)!} \sum_{m=0}^{n-i} \binom{n-i}{m} \frac{(-1)^m}{i+m} f_{_{T2GWG}}(x;(i+m)\alpha,\beta)
\end{align}
Therefore, we can present $f_{i:n}(x)$ as a linear combination of the T2GWG with parameter $(\alpha^*,\beta)$, where $\alpha^*=(i+m)\alpha$. 

\subsection{Stochastic ordering}
Stochastic ordering is a mathematical concept frequently applied in the realm of statistics, probability theory, decision theory, and economics \cite{Szekli2012Stochastic}. The most basic form of stochastic ordering is the usual order of real numbers, which extends naturally to random variables: a random variable $X$ is said to be stochastically smaller than another random variable $Y$ if, for every real number $x$, the probability that $X$ is less than or equal to $x$ is higher than or equal to the probability that $Y$ is less than or equal to $x$. This gives rise to the concept of one distribution being "stochastically larger" than another, which can be a valuable tool in comparing different probability models or assessing risk.
There are several types of stochastic orderings, such as increasing convex order, likelihood ratio order, and hazard rate order, each imposing a different structure on the sets of random variables or distributions.
Stochastic ordering is a significant concept because it enables us to make statements about the relative behavior of different random variables or distributions without specifying them precisely. It has been widely used in various fields, such as reliability, insurance, finance, operations research, and queueing theory. 

Let $X_1 \sim T2GWG(x;\alpha_1, \beta, \bpsi)$ and $X_2\sim T2GWG(x;\alpha_2, \beta, \bpsi)$. The likelihood ratio is 
\begin{align}
    \frac{f_{X_1}(x)}{f_{X_2}(x)} & = \frac{\alpha_1 \beta h(x,\bpsi) \frac{H(x,\bpsi)^{-\beta -1}}{\overline{H}(x,\bpsi)^{-\beta +1 }} \exp\left\{- \alpha_1 \left[ \frac{H(x,\bpsi)}{\overline{H}(x,\bpsi)} \right]^{-\beta} \right\}}{\alpha_2 \beta h(x,\bpsi) \frac{H(x,\bpsi)^{-\beta -1}}{\overline{H}(x,\bpsi)^{-\beta +1}} \exp\left\{- \alpha_2  \left[\frac{H(x,\bpsi)}{\overline{H}(x,\bpsi)} \right]^{-\beta}\right\}} \nonumber \\
    & = \frac{\alpha_1}{\alpha_2} \exp \left\{(\alpha_2-\alpha_1) \left[ \frac{H(x,\bpsi)}{\overline{H}(x,\bpsi)}\right]^{-\beta} \right\}
    \label{lr_alpha}
\end{align}
Then we differentiate the Eqn.\ref{lr_alpha} and obtain
\begin{align}
    \frac{d}{dx}\left(\frac{f_{X_1}(x)}{f_{X_2}(x)}\right) & = \frac{\alpha_1}{\alpha_2} \exp \left\{(\alpha_2-\alpha_1) \left[ \frac{H(x,\bpsi)}{\overline{H}(x,\bpsi)} \right]^{-\beta} \right\} \nonumber \\
    & \times \left[(\alpha_2-\alpha_1) h(x,\bpsi) \frac{H(x,\bpsi)^{-\beta-1}}{\overline{H}(x,\bpsi)^{-\beta+1}} \right]
\end{align}
If $\alpha_1 < \alpha_2$, $\frac{d}{dx}\left( \frac{f_{X_1}(x)}{f_{X_2}(x)}\right) < 0$. Thus it indicates that $X \preceq_{lr} Y$. According to the theorem proposed by \cite{Szekli2012Stochastic}, both $X \preceq_{hr} Y$ and $X \preceq Y$ hold.

\section{Methods of Estimation}
\subsection{Maximum Likelihood Estimation}
We can estimate the unknown parameters of the Type-2 Gumble Eibull-G family distributions by using the widely used Maximum Likelihood Estimation (MLE). Let $\mathbf{\Delta}=(\alpha,\beta, \bpsi)^T$. Then the log-likelihood for $\mathbf{\Delta}$ is defined by
\begin{align}
    \ell(\mathbf{\Delta})&= n\log (\alpha)+n\log(\beta)+\sum_{i=1}^n\log h(x_i; \bpsi) -(\beta+1) \sum_{i=1}^n \log H(x_i,\bpsi) \nonumber \\
    & + (\beta-1)\sum_{i=1}^n \log[1- H(x_i,\bpsi)] -\alpha \sum_{i=1}^n \left[ \frac{H(x_i,\bpsi)}{1-H(x_i,\bpsi) }\right]^{-\beta}
\end{align}
The first derivative of $\ell(\mathbf{\Delta})$ with respect to $\mathbf{\Delta}$ are shown as following
\begin{align}
    \frac{\partial\ell}{\partial \alpha}=\frac{n}{\alpha}-\sum_{i=1}^n \left[\frac{H(x_i, \bpsi)}{1-H(x_i,\bpsi)}\right]^{-\beta}
\end{align}
\begin{align}
    \frac{\partial \ell}{\partial \beta} =&\frac{n}{\beta}- \sum_{i=1}^n\log H(x_i,\bpsi) + \sum_{i=1}^n\log [1-H(x_i,\bpsi)] \nonumber \\
    & +\alpha \sum_{i=1}^n \left[ \frac{H(x_i,\bpsi)}{1-H(x_i,\bpsi)}\right]^{-\beta} \log \left (\frac{H(x_i,\bpsi)}{1-H(x_i,\bpsi)} \right)
\end{align}
and
\begin{align}
    \frac{\partial \ell}{\partial \psi_k} =&\sum_{i=1}^n \frac{1}{h(x_i,\bpsi)} \frac{\partial h(x_i,\bpsi)}{\partial \psi_k} - (\beta+1) \sum_{i=1}^n \frac{1}{H(x_i,\bpsi)} \frac{\partial H(x_i,\bpsi)}{\partial\psi_k} \nonumber \\
    &-(\beta-1)\sum_{i=1}^n \frac{1}{1-H(x_i,\bpsi)} \frac{\partial H(x_i,\bpsi)}{\partial\psi_k} \nonumber \\
    &+\alpha\beta\sum_{i=1}^n \frac{H(x_i,\bpsi)^{-\beta-1}}{[1-H(x_i,\bpsi)]^{-\beta+1}} \frac{\partial H(x_i,\bpsi)}{\partial \psi_k} 
\end{align}
where $\psi_k$ is the $k^{th}$ element of the vector $\bpsi$.\\
We can maximize the log-likehood function $\ell(\mathbf{\Delta})$ by solving the nonlinear equations $\left( \frac{\partial \ell}{\partial \alpha}, \frac{\partial \ell}{\partial \beta}, \frac{\partial \ell}{\partial \psi_k} \right) = 0$ with numerical methods such as Newton–Raphson approach.

\subsection{Least Square and Weighted Least Square Estimation}
The least squares (LS) method is a commonly used technique in regression analysis for approximating the solution of overdetermined systems. The method provides the best linear unbiased estimates of the unknown parameters if the errors are homoscedastic and uncorrelated. On the other hand, the weighted least squares (WLSE) approach extends the least squares technique by incorporating the different variances of the observations. This method assigns a weight to each data point based on the variance of its error term, placing less emphasis on the observations with higher variances to make the overall model more reliable.

The LSE and WLSE techniques can also provide estimators in the model. 
The LS estimation is given by
\begin{align}
    LS(\mathbf{\Delta})&= \sum_{i=1}^n \left(F(x_i,\mathbf{\Delta})-\frac{i}{n+1}\right)^2 \nonumber \\
    & = \sum_{i=1}^n \left(\exp\left\{-\alpha \left[ \frac{H(x_i,\bpsi)}{\overline{H}(x_i,\bpsi)} \right]^{-\beta} \right\}-\frac{i}{n+1}\right)^2  \label{LSE}
\end{align}
By differentiating Eqn.\ref{LSE}, we have the followings
\begin{align}
    \frac{\partial LS}{\partial \alpha} =& 2\sum_{i=1}^n \left(\exp\left\{-\alpha \left[ \frac{H(x_i,\bpsi)}{\overline{H}(x_i,\bpsi)} \right]^{-\beta} \right\}-\frac{i}{n+1}\right) \nonumber \\
    & \times (-1) \left[ \frac{H(x_i,\bpsi)}{\overline{H}(x_i,\bpsi)} \right]^{-\beta} \exp\left\{-\alpha \left[ \frac{H(x_i,\bpsi)}{\overline{H}(x_i,\bpsi)} \right]^{-\beta} \right\}
\end{align}
\begin{align}
    \frac{\partial LS}{\partial \beta} = & 2\sum_{i=1}^n \left(\exp\left\{-\alpha \left[\frac{H(x_i,\bpsi)}{\overline{H}(x_i,\bpsi)} \right]^{-\beta} \right\}-\frac{i}{n+1}\right) \nonumber \\
    & \times \alpha \log \left(\frac{H(x_i,\bpsi)}{\overline{H}(x_i,\bpsi)}\right) \left[ \frac{H(x_i,\bpsi)}{\overline{H}(x_i,\bpsi)} \right]^{-\beta} \exp\left\{-\alpha \left[ \frac{H(x_i,\bpsi)}{\overline{H}(x_i,\bpsi)} \right]^{-\beta} \right\}
\end{align}
and
\begin{align}
    \frac{\partial LS}{\partial \psi_k} =& 2\sum_{i=1}^n \left(\exp\left\{-\alpha \left[ \frac{H(x_i,\bpsi)}{\overline{H}(x_i,\bpsi)} \right]^{-\beta} \right\}-\frac{i}{n+1}\right) \nonumber \\
    & \times \alpha\beta\frac{\partial H(x_i,\bpsi)}{\partial \psi_k} \frac{H(x_i,\bpsi)^{-\beta-1}}{\overline{H}(x_i,\bpsi)^{-\beta +1}} \exp\left\{-\alpha \left[ \frac{H(x_i,\bpsi)}{\overline{H}(x_i,\bpsi)}\right]^{-\beta} \right\}
\end{align}
We can also apply Newton–Raphson procedure to minimize the least square estimation $LS(\mathbf{\Delta})$ by solving the equations $\left( \frac{\partial LS}{\partial \alpha}, \frac{\partial LS}{\partial \beta}, \frac{\partial LS}{\partial \theta_k} \right) = 0$ \\
Similarly, the WLS estimation can be obtained by minimizing
\begin{align}
    WLS(\mathbf{\Delta})& = \sum_{i=1}^n \frac{(n+1)^2(n+2)}{i(n-i+1)} \left(F(x_i, \mathbf{\Delta})-\frac{i}{n+1}\right)^2   \nonumber \\
    & = \sum_{i=1}^n \frac{(n+1)^2(n+2)}{i(n-i+1)} \left(\exp\left\{-\alpha \left[ \frac{H(x_i,\bpsi)}{\overline{H}(x_i,\bpsi)} \right]^{-\beta} \right\}-\frac{i}{n+1}\right)^2
\end{align}

\subsection{Maximum Product Spacing Approach of Estimation}

The Maximum Product Spacing (MPS) approach is particularly useful when dealing with unknown or complex distributions \cite{ChengAminMPS}. 
Unlike MLE, the MPS method does not require the explicit formulation of a likelihood function, making it a versatile and robust approach for different types of distributions.
The geometric mean of the MPS spacings is given by
\begin{align}
    G(\mathbf{\Delta}) = \left\{\prod_{i=1}^{n+1} D_i(x_i,\mathbf{\Delta}) \right\}^{\frac{1}{n+1}}
\end{align}
where 
\begin{align}
    D_i = \left\{
    \begin{aligned}
    \nonumber
    & F(x_1,\mathbf{\Delta}),  i=0\\
    & F(x_i,\mathbf{\Delta})-F(x_{i-1},\mathbf{\Delta}), i=2,3,...,n \\
    & 1-F(x_{n},\mathbf{\Delta}), i=n+1
    \end{aligned}
    \right.
\end{align}
Thus, We can maximize
\begin{align}
    G(\mathbf{\Delta}) =& \left[ \exp\left\{-\alpha \left[ \frac{H(x_1,\bpsi)}{\overline{H}(x_1,\bpsi)}\right]^{\beta}\right\} \left(1-\exp\left \{-\alpha\left[ \frac{H(x_n,\bpsi)}{\overline{H}(x_n,\bpsi)}\right]^{-\beta} \right\}\right)\right.\nonumber \\ 
    & \left. \times \prod_{i=2}^{n} \left(\exp\left\{-\alpha \left[ \frac{H(x_i,\bpsi)}{\overline{H}(x_i,\bpsi)} \right]^{-\beta} \right\}-\exp\left\{-\alpha\left[ \frac{H(x_{i-1},\bpsi)}{\overline{H}(x_{i-1},\bpsi)} \right]^{-\beta} \right\}\right) \right]^{\frac{1}{n+1}}
\end{align}
Equivalently, we can also maximize $H = \log G$
\begin{align}
    H(\mathbf{\Delta}) =& \frac{1}{n+1} \sum_{i=1}^{n+1} \log D_i(x_i,\mathbf{\Delta})\nonumber\\
    =& \frac{1}{n+1} \left\{-\alpha \left[\frac{H(x_1,\bpsi)}{\overline{H}(x_1,\bpsi)} \right]^{-\beta} + \log \left(1-\exp\left\{-\alpha\left[ \frac{H(x_n,\bpsi)}{\overline{H}(x_n,\bpsi)}\right]^{-\beta} \right\}\right)\right.\nonumber \\
    & \left. + \sum_{i=2}^n \log \left(\exp\left\{-\alpha \left[ \frac{H(x_i,\bpsi)}{\overline{H}(x_i,\bpsi)} \right]^{-\beta} \right\}-\exp\left\{-\alpha \left[ \frac{H(x_{i-1},\bpsi)}{\overline{H}(x_{i-1},\bpsi)} \right]^{-\beta} \right\}\right) \right\}
\end{align}
By solving $\left( \frac{\partial H}{\partial \alpha}, \frac{\partial H}{\partial \beta}, \frac{\partial H}{\partial \psi_k} \right) = 0$, we can obtain the estimators. The partial derivatives are provided in Appendix-\ref{sec:MPS_pd}.
\subsection{Cram\'er-von Mises Approach of Estimation}
The Cram\'er–von Mises method is another approach to estimate the parameters of a distribution \cite{MacDonaldCCVM}. The Cramér-von-Mises statistic measures the difference between the empirical distribution function of the data and the cumulative distribution function of the proposed model. This technique has an advantage over methods such as maximum likelihood estimation in that it considers the whole data set, not just the location and dispersion, resulting in a more comprehensive estimation.
We can apply the Cram\'er–von Mises criterion to obtain the estimators by minimizing the function $S(x;\mathbf{\Delta})$ with respect to $\mathbf{\Delta}$, where 
\begin{align}
    S(x;\mathbf{\Delta})&=\frac{1}{12n^2}+\frac{1}{n}\sum_{i=1}^n \left(F(x_i; \mathbf{\Delta})-\frac{2i-1}{2n} \right)^2 \nonumber \\
    & = \frac{1}{12n^2}+\frac{1}{n}\sum_{i=1}^n \left(\exp\left\{-\alpha \left[ \frac{H(x_i,\bpsi)}{\overline{H}(x_i,\bpsi)} \right]^{-\beta} \right\}-\frac{2i-1}{2n}\right)^2
\end{align}
Take the first partial derivatives of S, and we can have the following
\begin{align}
    \frac{\partial S}{\partial \alpha} =& \frac{2}{n}\sum_{i=1}^n \left(\exp\left\{-\alpha \left[ \frac{H(x_i,\bpsi)}{\overline{H}(x_i,\bpsi)}\right]^{-\beta} \right\}-\frac{2i-1}{2n}\right) \nonumber \\
    & \times (-1) \left[ \frac{H(x_i,\bpsi)}{\overline{H}(x_i,\bpsi)} \right]^{-\beta} \exp\left\{-\alpha \left[ \frac{H(x_i,\bpsi)}{\overline{H}(x_i,\bpsi)} \right]^{-\beta} \right\}
\end{align}
\begin{align}
    \frac{\partial S}{\partial \beta}=& \frac{2}{n}\sum_{i=1}^n \left(\exp\left\{-\alpha \left[ \frac{H(x_i,\bpsi)}{\overline{H}(x_i,\bpsi)} \right]^{-\beta} \right\}-\frac{2i-1}{2n}\right) \nonumber \\
    & \times \alpha \log \left(\frac{H(x_i,\bpsi)}{\overline{H}(x_i,\bpsi)}\right) \left[ \frac{H(x_i,\bpsi)}{\overline{H}(x_i,\bpsi)} \right]^{-\beta} \exp\left\{-\alpha \left[ \frac{H(x_i,\bpsi)}{\overline{H}(x_i,\bpsi)} \right]^{-\beta} \right\}
\end{align}
\begin{align}
    \frac{\partial S}{\partial \psi_k}=& \frac{2}{n}\sum_{i=1}^n \left(\exp\left\{-\alpha \left[ \frac{H(x_i,\bpsi)}{\overline{H}(x_i,\bpsi)} \right]^{-\beta} \right\}-\frac{2i-1}{2n}\right) \nonumber \\
    & \times \alpha\beta\frac{\partial H(x_i,\bpsi)}{\partial \psi_k} \frac{H(x_i,\bpsi)^{-\beta -1}}{\overline{H}(x_i,\bpsi)^{-\beta +1}}\exp\left\{-\alpha \left[ \frac{H(x_i,\bpsi)}{\overline{H}(x_i,\bpsi)} \right]^{-\beta} \right\}
\end{align}

\subsection{Anderson and Darling Approach of Estimation}
The Anderson-Darling approach was proposed by \cite{anderson1954test} to test if a data set follows a specific distribution, and can also be used to estimate parameters.
This method gives greater weight to the tails of the distribution
compared to other methods, like the Kolmogorov-Smirnov test. The Anderson-Darling statistic is minimized to find the parameters of the best-fitting distribution. This estimation method is highly sensitive to deviations in the tails and thus can be more powerful for identifying whether a particular distribution fits the data.
The Anderson-Darling estimators can be obtained by minimizing
\begin{align}
    AD(\mathbf{\Delta})  =& -n-\frac{1}{n} \sum_{i=1}^n (2i-1)[\log F(x_i, \mathbf{\Delta})- \log \overline{F}(x_{n+1-i}, \mathbf{\Delta})] \nonumber \\
     = & -n-\frac{1}{n} \sum_{i=1}^n (2i-1) \left[-\alpha \left[ \frac{H(x_i,\bpsi)}{\overline{H}(x_i,\bpsi)} \right]^{-\beta} \right. \nonumber \\
    & - \left. \log \left(1-\exp\left\{-\alpha \left[ \frac{H(x_{n+1-i},\bpsi)}{\overline{H}(x_{n+1-i},\bpsi)} \right]^{-\beta} \right\} \right)\right]
\end{align}
Similarly, we take the first derivatives of $AD(\mathbf{\Delta})$ and obtain
\begin{align}
    \frac{\partial AD}{\partial \alpha} =& -\frac{1}{n} \sum_{i=1}^n (2i-1) \left(-\left[ \frac{H(x_i,\bpsi)}{\overline{H}(x_i,\bpsi)} \right]^{-\beta} \right. \nonumber \\
    & - \left. \left[ \frac{H(x_{n+1-i},\bpsi)}{\overline{H}(x_{n+1-i},\bpsi)} \right]^{-\beta} \frac{\exp\left\{-\alpha \left[ \frac{H(x_{n+1-i},\bpsi)}{\overline{H}(x_{n+1-i},\bpsi)} \right]^{-\beta} \right\}}{1-\exp\left\{-\alpha \left[\frac{H(x_{n+1-i},\bpsi)}{\overline{H}(x_{n+1-i},\bpsi)} \right]^{-\beta} \right\}} \right)
\end{align}
\begin{align}
    \frac{\partial AD}{\partial \beta} =& -\frac{1}{n} \sum_{i=1}^n (2i-1)\left[\alpha\log\left(\frac{H(x_i,\bpsi)}{\overline{H}(x_i,\bpsi)} \right)\left[\frac{H(x_i,\bpsi)}{\overline{H}(x_i,\bpsi)} \right]^{-\beta} \right.  \nonumber \\
    &+ \alpha \log \left( \frac{H(x_{n+1-i},\bpsi)}{\overline{H}(x_{n+1-i},\bpsi)} \right) \left[\frac{H(x_{n+1-i},\bpsi)}{\overline{H}(x_{n+1-i},\bpsi)} \right]^{-\beta} 
    \frac{\exp\left\{-\alpha \left[ \frac{H(x_{n+1-i},\bpsi)}{\overline{H}(x_{n+1-i},\bpsi)} \right]^{-\beta} \right\}}{1-\exp\left\{-\alpha \left[ \frac{H(x_{n+1-i},\bpsi)}{\overline{H}(x_{n+1-i},\bpsi)} \right]^{-\beta} \right\}} 
\end{align}
\begin{align}
    \frac{\partial AD}{\partial \psi_k}=&-\frac{1}{n}\sum_{i=1}^n (2i-1) \left[\alpha \beta \frac{H(x_i,\bpsi)^{-\beta-1}}{\overline{H}(x_i,\bpsi)^{-\beta +1}} \frac{\partial H(x_i,\bpsi)}{\partial\psi_k}\right.\nonumber \\
    & \left. +\alpha \beta \frac{H(x_{n+1-i},\bpsi)^{-\beta-1}}{\overline{H}(x_{n+1-i},\bpsi)^{-\beta+1}} \frac{\partial H(x_{n+1-i},\bpsi)}{\partial \psi_k}  \frac{\exp\left\{-\alpha \left[ \frac{H(x_{n+1-i},\bpsi)}{\overline{H}(x_{n+1-i},\bpsi)} \right]^{-\beta} \right\}}{1-\exp\left\{-\alpha \left[ \frac{H(x_{n+1-i},\bpsi)}{\overline{H}(x_{n+1-i},\bpsi)} \right]^{-\beta} \right\}}
    \right] 
\end{align}

\subsection{Simulation and Estimation}
\begin{figure}[htbp]
    \centering
    \begin{minipage}{0.48\linewidth}
	\centering
        \includegraphics[width=0.9\linewidth]{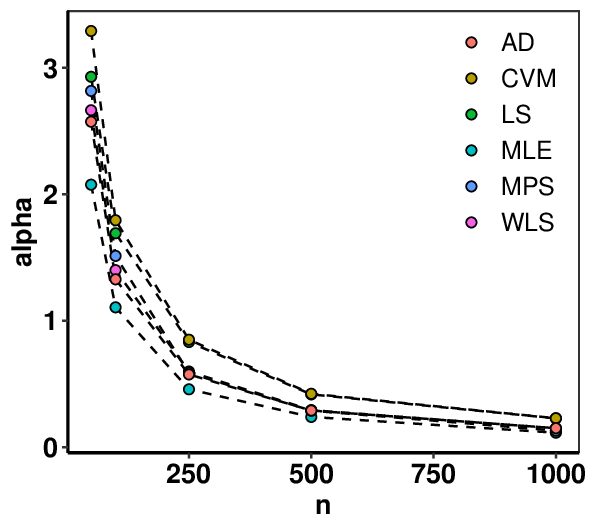}
	\label{alpha_est}
    \end{minipage}
    \begin{minipage}{0.48\linewidth}
	\centering
	\includegraphics[width=0.9\linewidth]{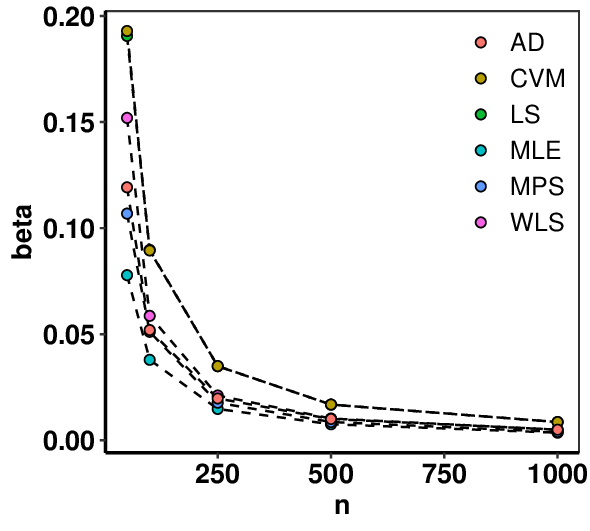}
	\label{beta_est}
    \end{minipage}
    \qquad

    \begin{minipage}{0.48\linewidth} \centering
       \includegraphics[width=0.9\linewidth]{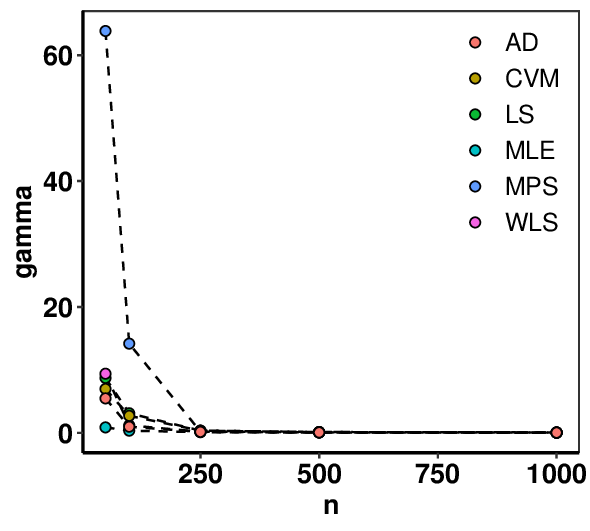}
       \label{gamma_est}
    \end{minipage}
    \caption{MSE of parameters in Table 3}
    \label{MSE_plot}
\end{figure}

We combined Monte Carlo simulation with the above techniques to estimate the parameters. The parameters are set as $\alpha=2.5, \beta=0.8$, and $\gamma=1.3$. The sample sizes $N=50, 100, 250, 500$, and $1000$ were used to generate random samples. For each sample size, the experiment was replicated for $N=1000$ times. Then the bias and mean squared error (MSE) were calculated. Table~\ref{est_results_new} and Fig.~\ref{MSE_plot} show the estimation results. The MSE converges to $0$ when $N$ increases, confirming the estimations' stability and reliability.

\begin{sidewaystable}[thp]
\caption{Different estimations for $\alpha=2.5,\beta=0.8,\gamma=1.3$}
\begin{tabular}{ccccccccccccccc}
\toprule
\multicolumn{2}{c}{} &\multicolumn{2}{c}{\textbf{MLE}} &\multicolumn{2}{c}{\textbf{LS}} &\multicolumn{2}{c}{\textbf{WLS}} &\multicolumn{2}{c}{\textbf{MPS}} &\multicolumn{2}{c}{\textbf{CVM}} &\multicolumn{2}{c}{\textbf{AD}} \\
\midrule
\textbf{$N$}&&\textbf{Bias} & \textbf{MSE} & \textbf{Bias} & \textbf{MSE}& \textbf{Bias} & \textbf{MSE}& \textbf{Bias} & \textbf{MSE}& \textbf{Bias} & \textbf{MSE}& \textbf{Bias} & \textbf{MSE} \\
\midrule
\multirow{3}{*}{50}&$\alpha$&0.0909	&	2.0761	&	-0.0418	&	2.9292	&	0.0307	&	2.6645	&	0.8349	&	2.8174	&	0.1095	&	3.2913	&	0.1459	&	2.5743 \\
&$\beta$&  0.0512	&	0.0778	&	0.0107	&	0.1906	&	0.0066	&	0.1520	&	-0.1492	&	0.1068	&	0.0262	&	0.1929	&	0.0010	&	0.1192 \\
&$\gamma$&  0.1512	&	0.8690	&	0.9373	&	8.7264	&	0.8477	&	9.4008	&	2.4409	&	63.8413	&	0.8518	&	6.9832	&	0.6359	&	5.4760\\
\midrule
\multirow{3}{*}{100}&$\alpha$&0.0662	&	1.1061	&	0.0595	&	1.6917	&	0.0878	&	1.4001	&	0.5914	&	1.5128	&	0.1357	&	1.7947	&	0.1327	&	1.3275 \\
&$\beta$&0.0182	&	0.0379	&	-0.0229	&	0.0897	&	-0.0158	&	0.0586	&	-0.1047	&	0.0511	&	-0.0161	&	0.0895	&	-0.0163	&	0.0520 \\
&$\gamma$&0.0963	&	0.3519	&	0.4901	&	3.1145	&	0.2684	&	1.1972	&	0.8138	&	14.1750	&	0.4755	&	2.7094	&	0.2518	&	0.9787 \\
\midrule
\multirow{3}{*}{250}&$\alpha$&0.0203	&	0.4578	&	0.0211	&	0.8317	&	0.0336	&	0.5983	&	0.3204	&	0.5872	&	0.0510	&	0.8502	&	0.0533	&	0.5754 \\
&$\beta$&0.0097	&	0.0148	&	-0.0027	&	0.0349	&	-0.0010	&	0.0211	&	-0.0544	&	0.0177	&	-0.0002	&	0.0350	&	-0.0027	&	0.0197 \\
&$\gamma$&0.0265	&	0.1085	&	0.1088	&	0.3554	&	0.0614	&	0.1707	&	0.2180	&	0.2020	&	0.1141	&	0.3529	&	0.0660	&	0.1620 \\
\midrule
\multirow{3}{*}{500}&$\alpha$&0.0226	&	0.2394	&	0.0168	&	0.4176	&	0.0220	&	0.2933	&	0.2066	&	0.2890	&	0.0318	&	0.4223	&	0.0321	&	0.2900 \\
&$\beta$&0.0049	&	0.0075	&	-0.0015	&	0.0168	&	0.0002	&	0.0102	&	-0.0335	&	0.0086	&	-0.0003	&	0.0168	&	-0.0008	&	0.0100 \\
&$\gamma$&0.0189	&	0.0543	&	0.0498	&	0.1280	&	0.0307	&	0.0742	&	0.1269	&	0.0809	&	0.0529	&	0.1282	&	0.0347	&	0.0744 \\
\midrule
\multirow{3}{*}{1000}&$\alpha$&0.0072	&	0.1146	&	0.0071	&	0.2278	&	0.0113	&	0.1497	&	0.1163	&	0.1310	&	0.0146	&	0.2290	&	0.0167	&	0.1504 \\
&$\beta$&0.0027	&	0.0036	&	-0.0003	&	0.0086	&	0.0001	&	0.0050	&	-0.0198	&	0.0039	&	0.0003	&	0.0086	&	-0.0007	&	0.0050 \\
&$\gamma$&0.0079	&	0.0254	&	0.0258	&	0.0673	&	0.0160	&	0.0373	&	0.0691	&	0.0331	&	0.0274	&	0.0673	&	0.0186	&	0.0380 \\
\bottomrule
\end{tabular}
\label{est_results_new}
\end{sidewaystable}

\section{Special Cases}
In this section, we will explore a variety of special cases that emerge from our novel distribution model. By closely examining these unique instances, we aim to illustrate the multifaceted aspects and potential applications of the T2GWG family of distributions. 

\subsection{Type-2 Gumbel Weibull-Exponential (T2GWE) distribution}
Suppose the baseline distribution $H(x,\bpsi)$ is an exponential distribution with parameter $\gamma > 0$. Then $h(x;\gamma) = \gamma e^{-\gamma x}$ and $H(x;\gamma)=1-e^{-\gamma x}$. 

\subsubsection{cdf and pdf of the T2GWE distribution}
The cdf of the Type-2 Gumbel Weibull-Exponential distribution is presented as 
\begin{align}
    F_{T2GWE}(x)=\exp\left\{- \alpha (e^{\gamma x}-1)^{-\beta} \right\}
\end{align}
and the pdf is given by
\begin{align}
    f_{T2GWE}(x)=\beta \alpha (e^{\gamma x}-1)^{-\beta-1}\exp \left\{ \gamma x-\alpha (e^{\gamma x}-1)^{-\beta} \right\}
\end{align}
where $x \geq 0$ and $\alpha,\beta> 0$

\subsubsection{Hazard rate and quantile functions}
\begin{figure}[htbp!]
    \centering
    \includegraphics[width=0.48\textwidth]{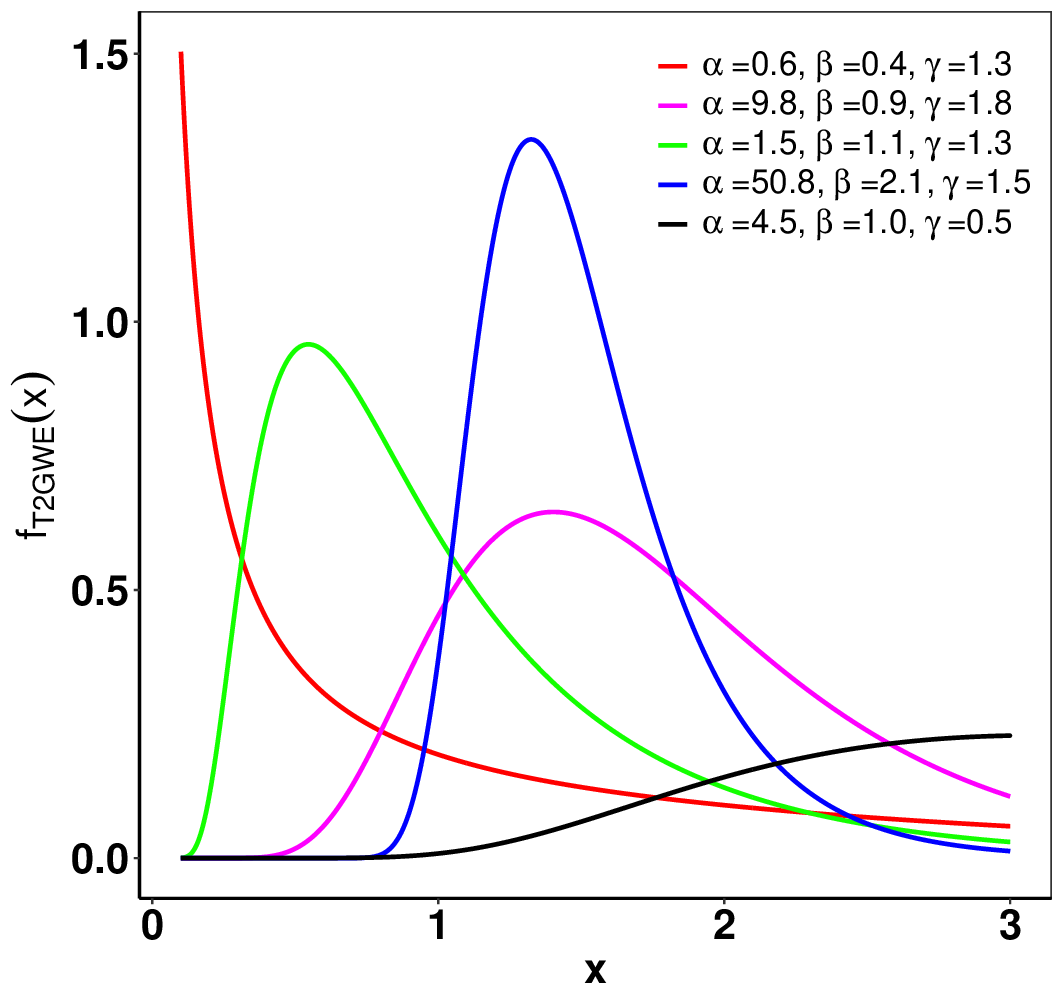}
    \includegraphics[width=0.48\textwidth]{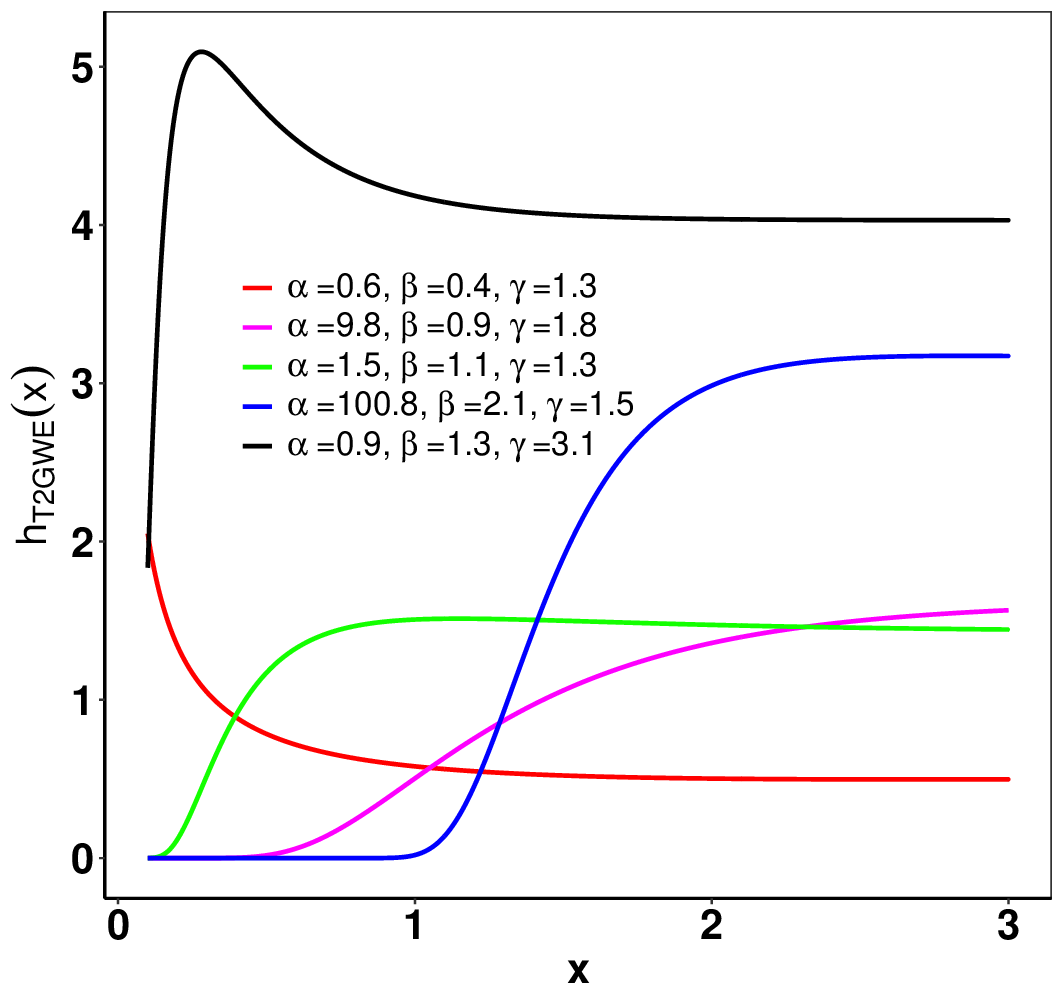}
    \caption{Left: The pdf of T2GWE distribution for different parameters. Right: The hrf of T2GWE for different parameter values.}
    \label{fig:t2gwe}
\end{figure}

The hrf is shown as
\begin{align}
    h_{T2GWE} = \frac{\beta \alpha (e^{\gamma x}-1)^{-\beta-1}\exp \left\{ \gamma x-\alpha (e^{\gamma x}-1)^{-\beta} \right\}}{1-\exp \left\{-\alpha(e^{\gamma x}-1)^{-\beta} \right\}}
\end{align}
and the reverse hrf is given by
\begin{align}
    \tau_{T2GWE} =\beta \alpha (e^{\gamma x}-1)^{-\beta-1}
\end{align}
Moreover, the quantile function is obtained as
\begin{align}
    x_p = -\frac{1}{\gamma}\log \left( \left[\frac{\log p}{-\alpha} \right]^{-\frac{1}{\beta}}+1\right)
\end{align}
Fig.~\ref{fig:t2gwe} displays several typical configurations of the pdf and hrf for the T2GWE distribution. The pdf of the T2GWE distribution shows various configurations, including almost symmetric, right-skewed, decreasing, and increasing. Additionally, the hrf of the T2GWE distribution can exhibit a range of shapes, such as decreasing, increasing, and right-skewed.

\subsection{Type-2 Gumbel Weibull-Uniform (T2GWU) distribution}
Let the baseline distribution $H(x,\bpsi)$ be a uniform distribution with parameter $\gamma > 0$. Then $h(x,\gamma) = \frac{1}{\gamma}$ and $H(x,\gamma) = \frac{x}{\gamma}$.
\subsubsection{cdf and pdf of the T2GWU distribution}
The cdf of the T2GWU distribution is 
\begin{align}
    F_{T2GWU}=\exp\left\{-\alpha\left( \frac{x}{\gamma-x}\right)^{-\beta} \right\},
\end{align}
and the corresponding pdf is 
\begin{align}
    f_{T2GWU}(x)=\alpha \beta \frac{\gamma x^{-\beta-1}}{(\gamma-x)^{-\beta+1}} \exp\left\{-\alpha\left( \frac{x}{\gamma-x}\right)^{-\beta}\right\}.
\end{align}
\subsubsection{Hazard rate and quantile functions}
The hrf of T2GWU is displayed by 
\begin{align}
    h_{T2GWU} = \frac{\alpha \beta \frac{\gamma x^{-\beta-1}}{(\gamma-x)^{-\beta+1}} \exp\left\{-\alpha \left(\frac{x}{\gamma-x}\right)^{-\beta}\right\}}{1-\exp \left\{- \alpha \left(\frac{x}{\gamma-x}\right)^{-\beta} \right\}},
\end{align}
and the reverse hrf is given by
\begin{align}
    \tau_{T2GWU} = \alpha \beta \frac{\gamma x^{-\beta-1}}{(\gamma-x)^{-\beta+1}}.
\end{align}
Moreover, the quantile function can be obtained as
\begin{align}
    x_p = \frac{1}{1+\left[\frac{\log p}{-\alpha}\right]^{\frac{1}{\beta}}}.
\end{align}

\begin{figure}[htbp!]
    \centering
    \includegraphics[width=0.48\textwidth]{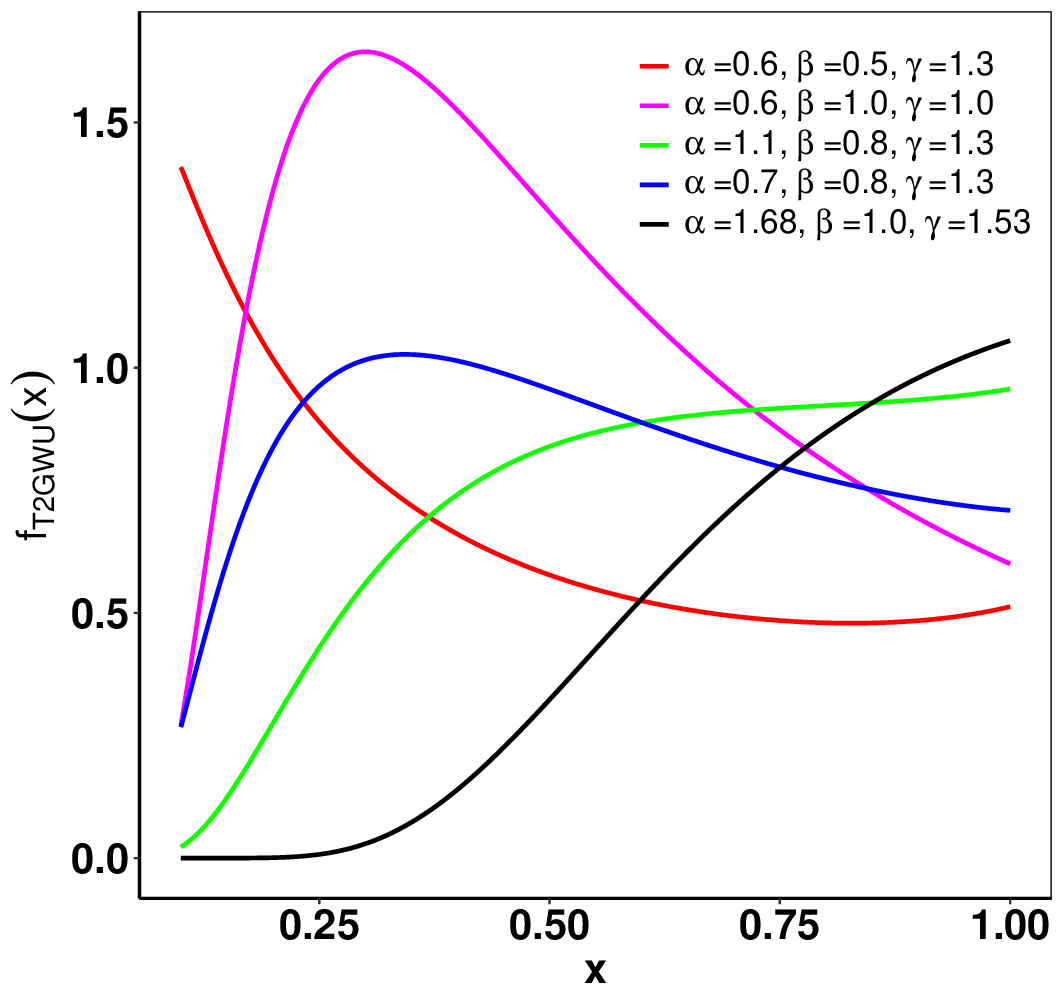}
     \includegraphics[width=0.48\textwidth]{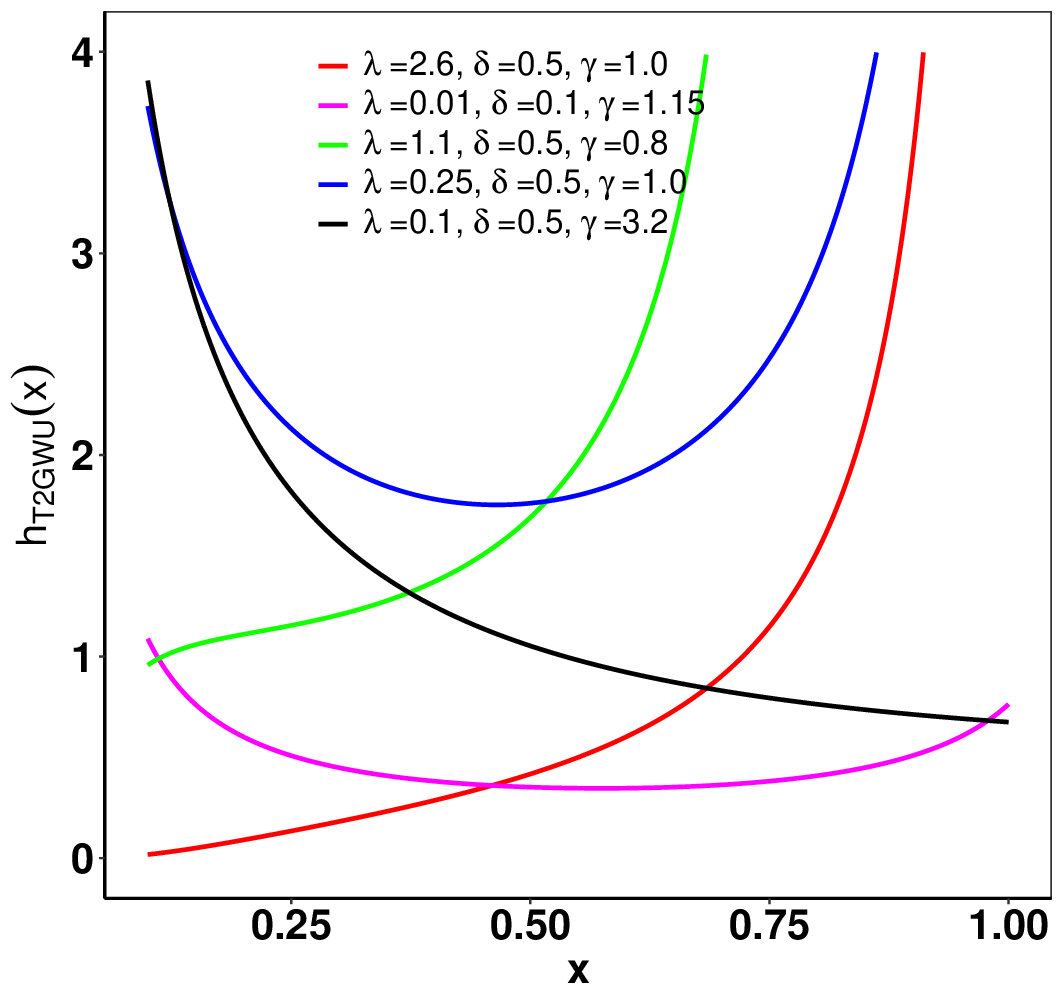}
    \caption{ Left: pdf of T2GWU distribution for different values of parameters $\alpha$, $\beta$,  and $\gamma$. Right: hrf  of T2GWU for selected parameters $\alpha$, $\beta$, and $\gamma$.}
    \label{fig:t2gwu}
\end{figure}

Fig.~\ref{fig:t2gwu} plots of the pdf and hrf for the T2GWU distribution with several combinations of parameter values. The pdf plots show different shapes, including right-skewed, decreasing, and increasing. In addition, the hrf plots capture various possibilities such as increasing, decreasing, bathtub, and shallow bathtub. 

\subsection{Type-2 Gumbel Weibull-Pareto (T2GWP) distribution}
If we set $H(x,\bpsi)$ as a Pareto distribution with parameter $\theta, k > 0$, then $\ds h(x,\theta, k) = \frac{k\theta^k}{x^{k+1}}$ and $\ds H(x,\theta, k) = 1-\left(\frac{\theta}{x}\right)^k$.
\subsubsection{cdf and pdf of the T2GWP distribution}
Thus the cdf of the T2GWP distribution is given by
\begin{align}
    F_{T2GWP} = \exp\left\{-\alpha \left(\left[\frac{x}{\theta} \right]^k-1\right)^{-\beta}\right\},
\end{align}
with a pdf
\begin{align}
    f_{T2GWP}(x)=\beta \alpha kx^{k-1}\theta^{-k} \left(\left[\frac{x}{\theta} \right]^k-1\right)^{-\beta-1}\exp\left\{-\alpha \left( \left[\frac{x}{\theta}\right]^k-1\right)^{-\beta}\right\}.
\end{align}

\subsubsection{Hazard rate and quantile functions}
The hrf of T2GWP is displayed by 
\begin{align}
    h_{T2GWP} = \frac{\beta \alpha kx^{k-1}\theta^{-k} \left(\left[ \frac{x}{\theta} \right]^k-1\right)^{-\beta-1}\exp\left\{-\alpha\left(\left[\frac{x}{\theta} \right]^k-1\right)^{-\beta} \right\}}{1-\exp\left\{-\alpha \left(\left[\frac{x}{\theta} \right]^k-1\right)^{-\beta} \right\}},
\end{align}
and the reverse hrf is
\begin{align}
    \tau_{T2GWP} = \beta \delta \alpha^{-\delta}kx^{k-1}\theta^{-k} \left(\left[\frac{x}{\theta} \right]^k-1\right)^{-\beta\delta-1}.
\end{align}
Moreover, the quantile function is obtained as
\begin{align}
    x_p=\theta\left(\left[\frac{\log p}{-\alpha}\right]^{-\frac{1}{\beta}} +1\right)^{\frac{1}{k}}
\end{align}

\begin{figure}[htbp!]
    \centering
    \includegraphics[width=0.48\textwidth]{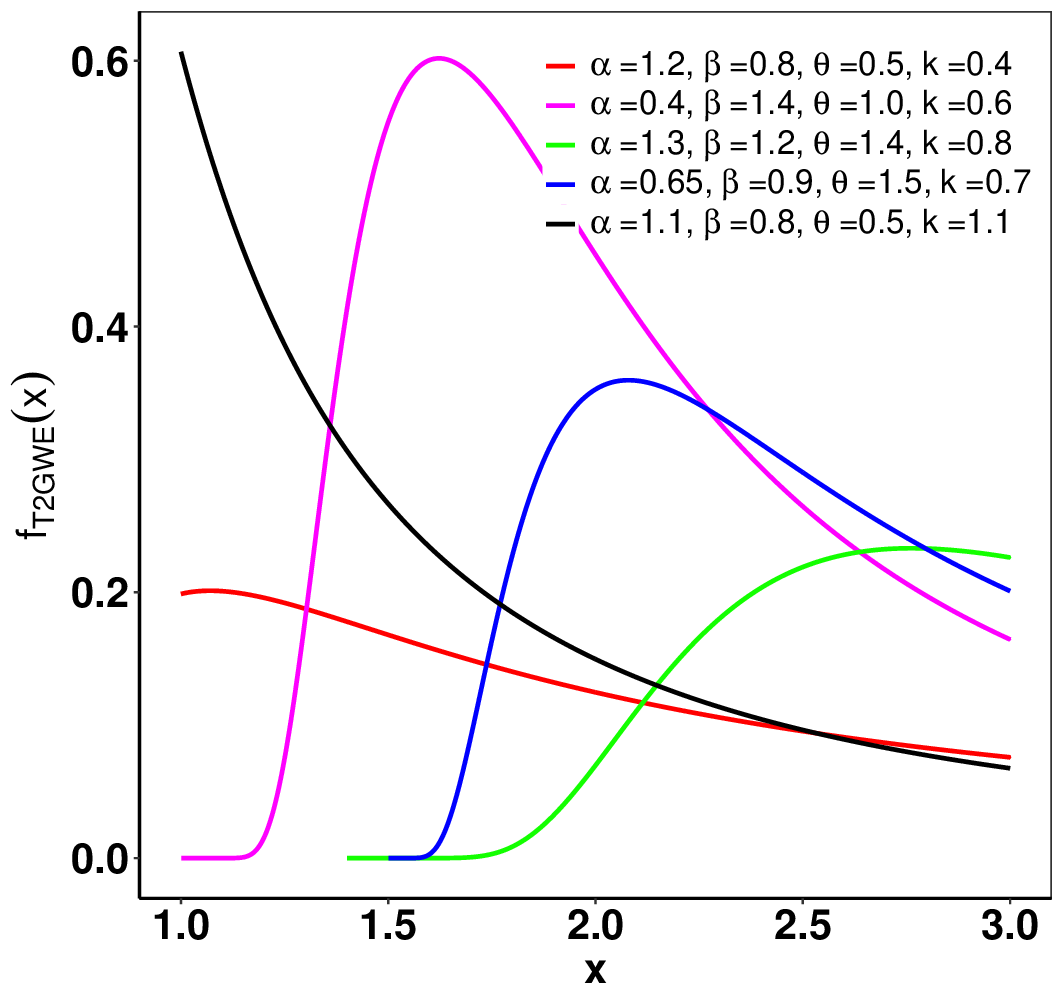}
    \includegraphics[width=0.48\textwidth]{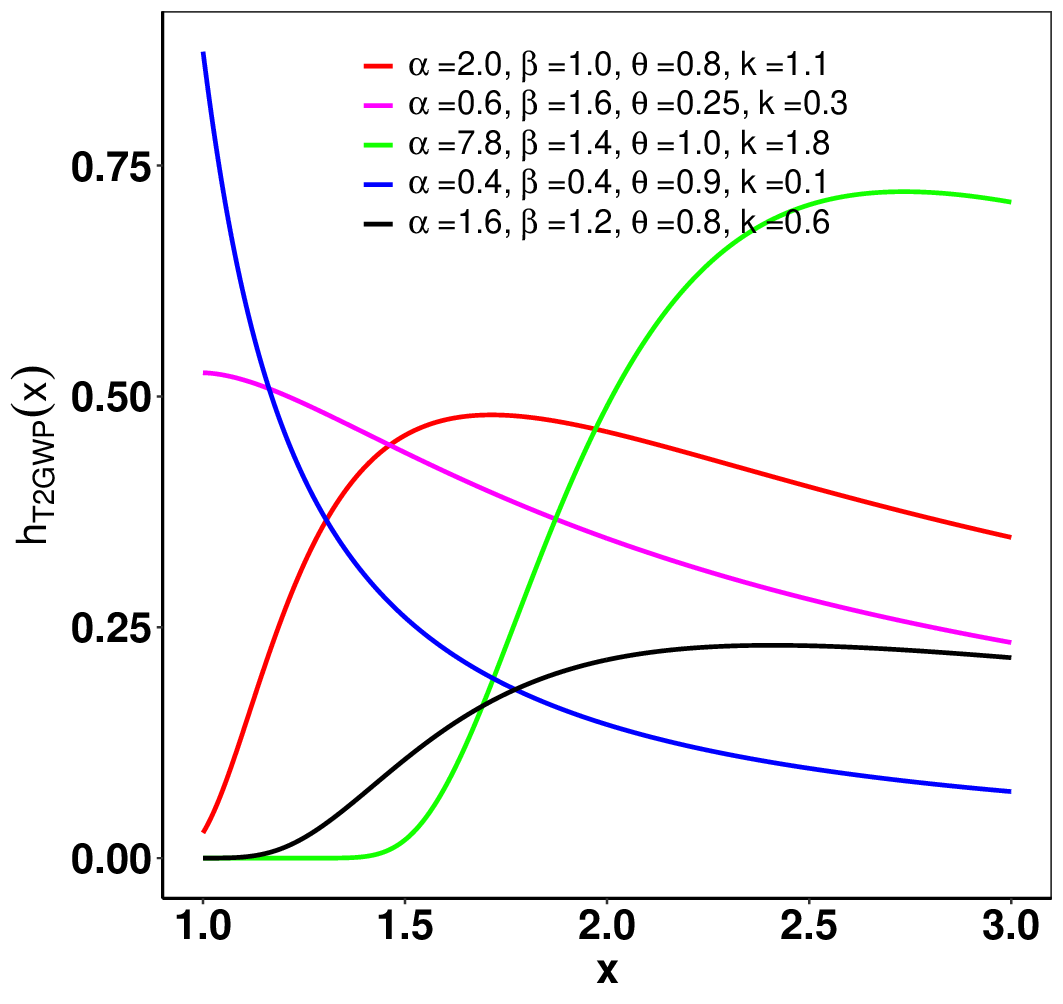}
    \caption{ Left: The pdf of T2GWP distribution for selected values of $\alpha$, $\beta$, $\theta$ and $k$. Right: The hrf of T2GWP for various $\alpha$, $\beta$, $\theta$ and $k$.}
    \label{fig:t2gwp}
\end{figure}

Shapes of the pdf and hrf for the T2GWP distribution with selected parameters are shown in Fig.~\ref{fig:t2gwp}. The pdfs exhibit a variety of shapes including right-skewed, decreasing, and increasing. Moreover, hrf plots for the T2GWP distribution display growing, decreasing, and right-skewed forms.

\section{Applications}

In this section, we will investigate beyond theoretical constructs and delve into the practical implications of our model, demonstrating its applicability using real-world data sets.
This will validate the practical utility of our newly devised model and shed light on how it can be effectively employed in handling concrete data-driven scenarios. 
The objective is to ensure that our theoretical advancements resonate with tangible applications, thereby significantly contributing to both the academic discourse and the operational applications of statistical distributions.

In this section, we present three applications of the Type-2 Gumbel Weibull-exponential distribution. We compared it with the Exponentiated Gumbel Type-2 (EGT) \cite{okorie2016exponentiated}, Weibull Generalized Exponential (WGE) \cite{mustafa2016weibull}, Lomax Gumbel Type-2 (LGT) \cite{adeyemi2022lomax}, Type-2 Gumbel (T2G), and with Exponentiated Weibull-Logistic distributions (EWL) \cite{murat2020exponentiated}. The pdf and cdf of those distributions are provided in Appendix~\ref{app:distributions}. The goodness-of-fit statistics including -2log-likelihood statistic, Cram\'er-von Mises statistic ($W^{*}$), Anderson-Darling statistic ($A^{*}$), Akaike Information Criterion (AIC), Bayesian Information Criterion (BIC), Consistent Akaike Information Criterion (CAIC), Hannan-Quinn criterion (HQIC), Kolmogorov-Smirnov test statistic (K-S) and its corresponding p-value are reported.  

To evaluate and compare the performance of different models, we can examine their goodness-of-fit statistics. Generally, a model with smaller values in these statistics fits the data better. However, it's important to note that for the p-values, which are a measure of expectation, a larger value indicates a better fit.

\subsection{Chemotherapy data}
This dataset is a subset of data reported by Bekker et al.  \cite{Bekker2000Chemo}, which represents the survival times (in years) of a group of patients who received chemotherapy treatment alone. Please refer to the ``Declarations" section for guidance on how to access the data.

\begin{figure}[htbp!]
    \centering
    \begin{minipage}{0.45\linewidth}
    \centering
    \includegraphics[width=0.9\linewidth]{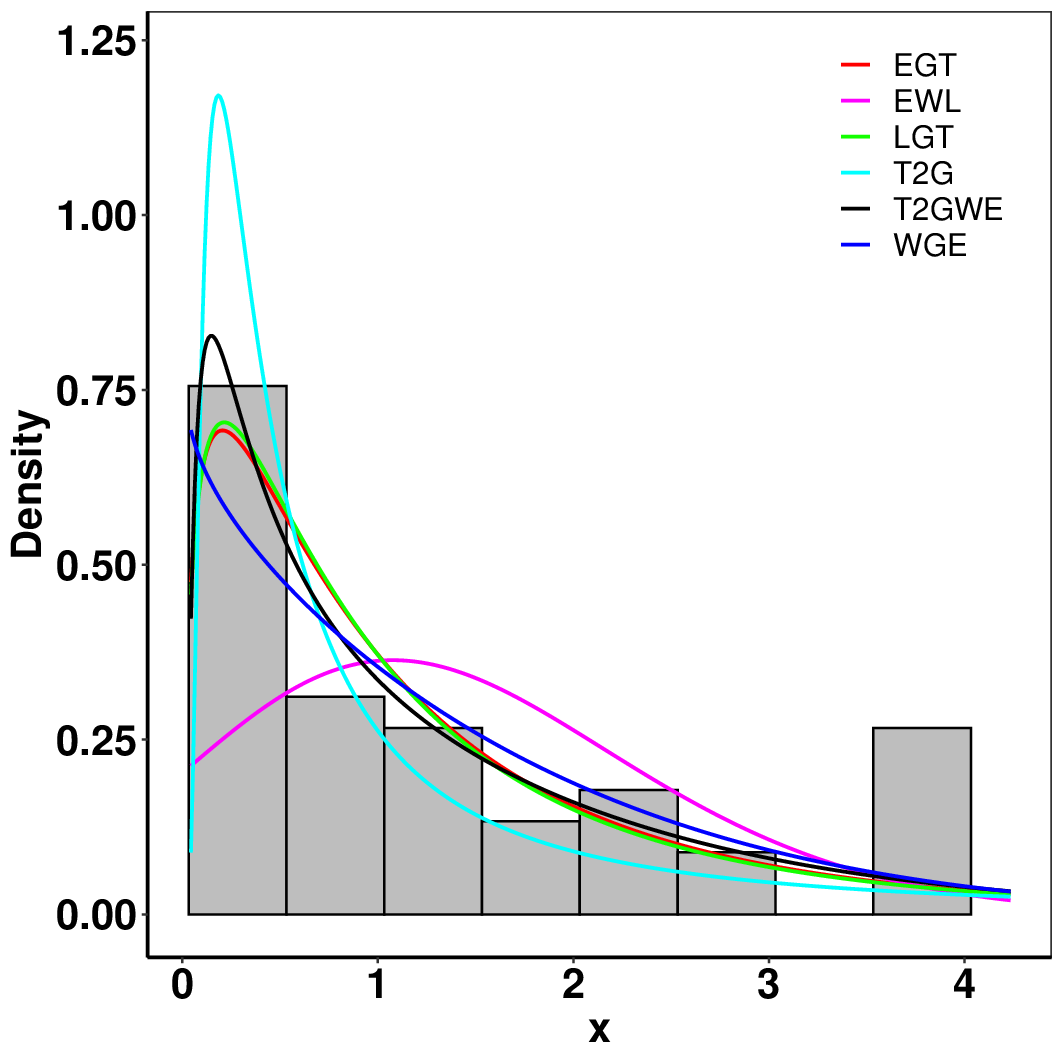}
    \end{minipage}
    \begin{minipage}{0.45\linewidth}
    \centering
    \includegraphics[width=0.9\linewidth]{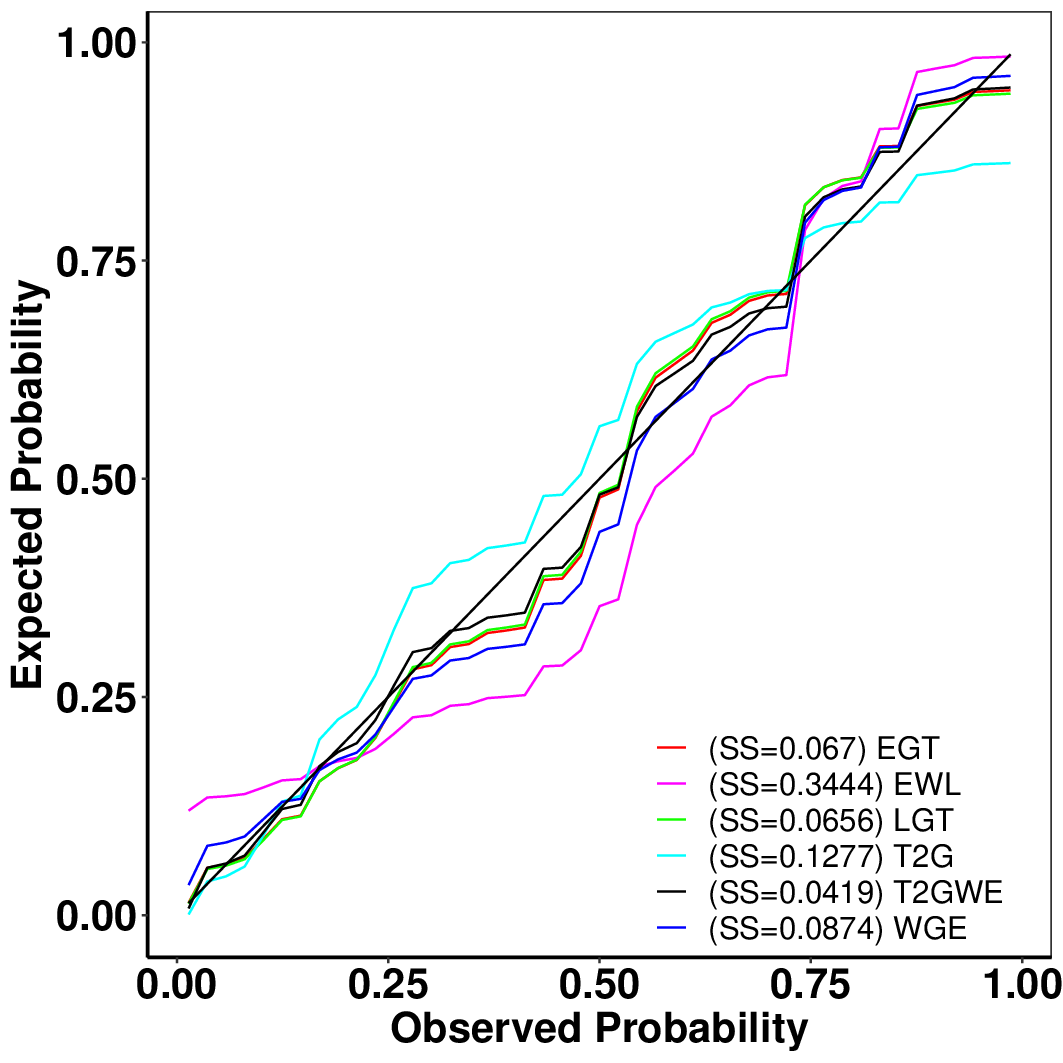}
    \end{minipage}
    \caption{left) Fitted density superposed on the histogram and observed probability for the Chemotherapy data. right) Expected probability plots for the Chemotherapy data.}
    \label{hist_ep_chemo}
\end{figure}

The estimates of the parameters and the goodness-of-fit statistics are summarized in Table~\ref{mle_compare_chemo}. Fig.~\ref{hist_ep_chemo} plots of the fitted densities alongside the histogram and the expected probability. From Table \ref{mle_compare_chemo}, we can conclude that the Type-2 Gumbel Weibull-exponential distribution has better performance than other distributions because of its lowest values among all measures of goodness-of-fit and the highest p-value in the K-S test. 

\begin{figure}[htbp!]
    \centering
    \begin{minipage}{0.45\linewidth}
    \centering
    \includegraphics[width=0.9\linewidth]{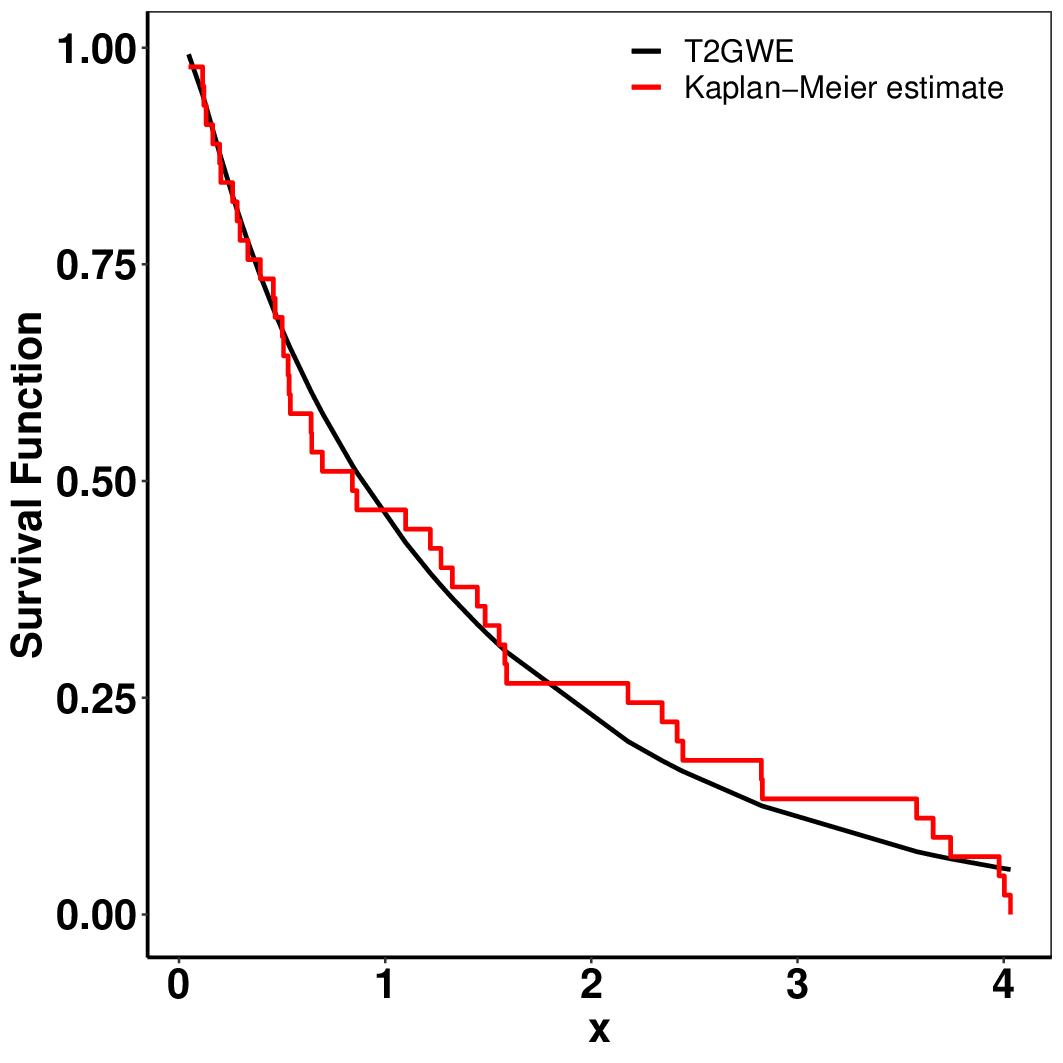}
    \end{minipage}
    \begin{minipage}{0.45\linewidth}
    \centering
    \includegraphics[width=0.9\linewidth]{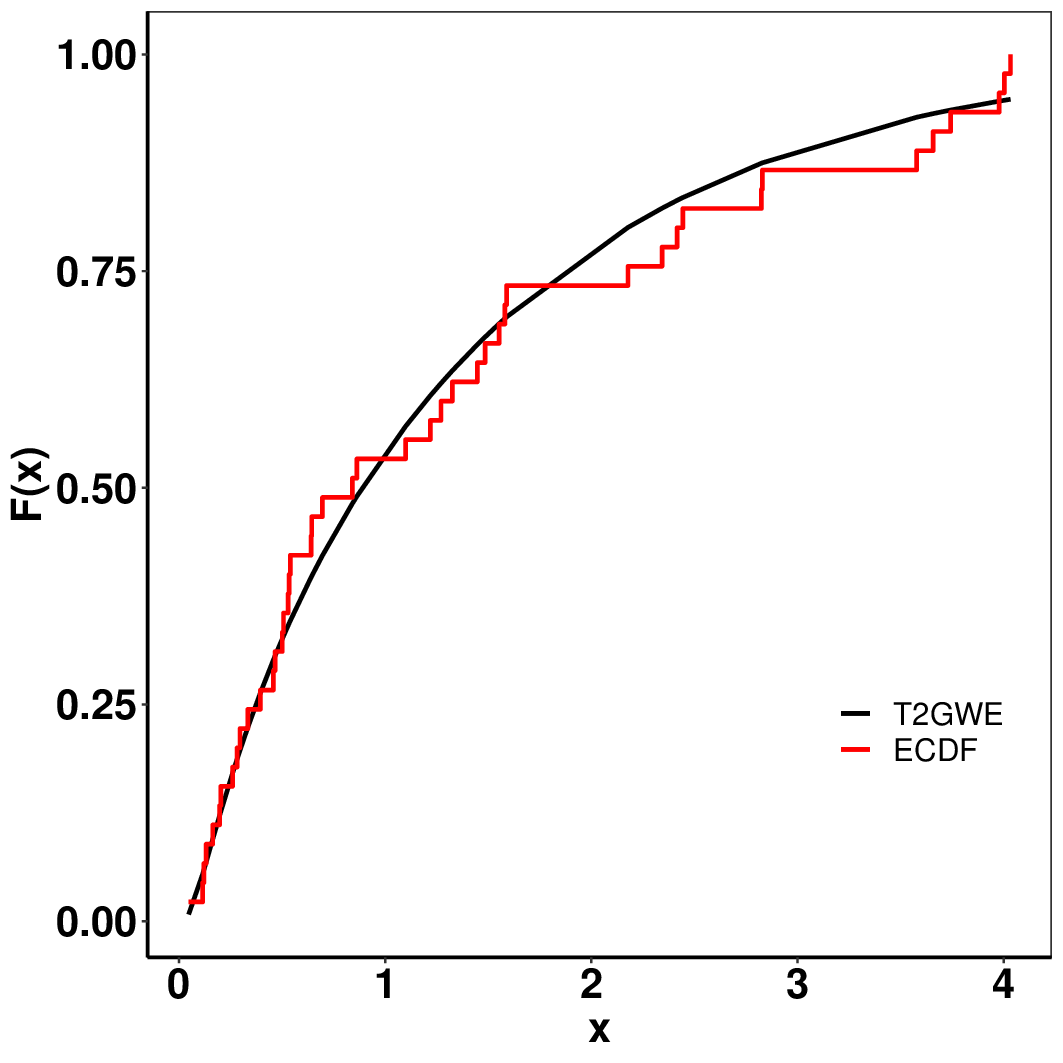}
    \end{minipage}
    \qquad
    
    \begin{minipage}{0.45\linewidth}
    \centering
    \includegraphics[width=0.9\linewidth]{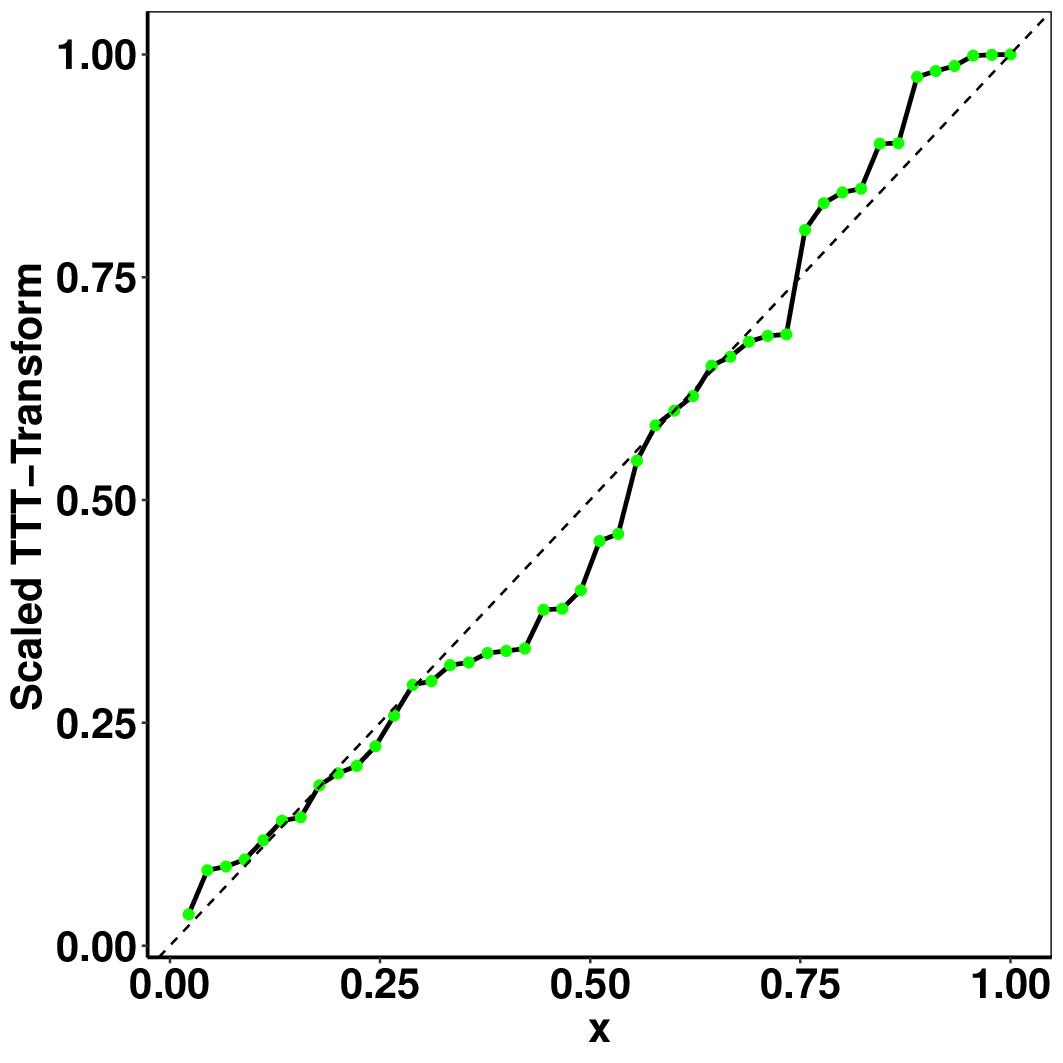}
    \end{minipage}
    \begin{minipage}{0.45\linewidth}
    \centering
    \includegraphics[width=0.9\linewidth]{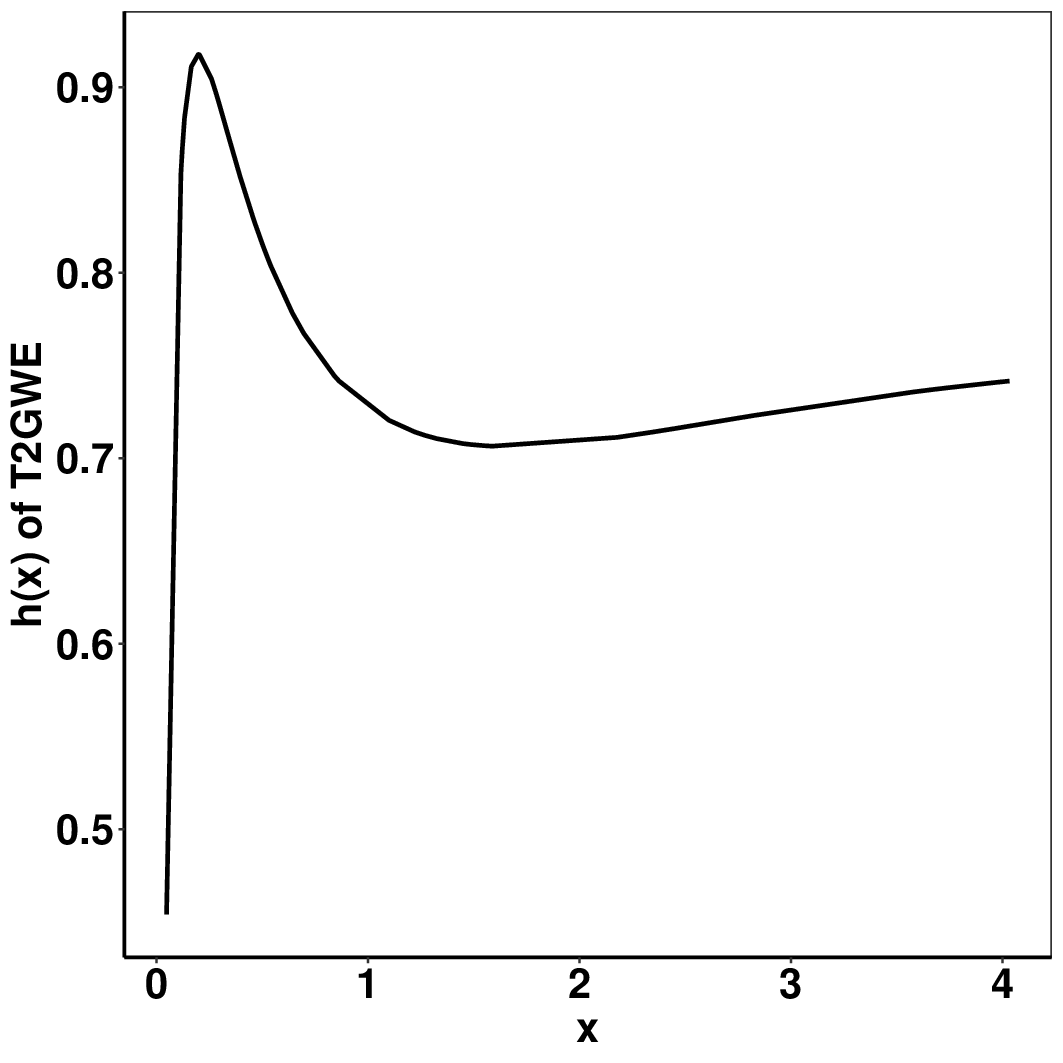}
    \end{minipage}
    \caption{Fitted K-M survival curve, theoretical and ECDF, the TTT statistics, and the hrf for the Chemotherapy data.}
    \label{chemo_4plot}
\end{figure}

In Fig.~\ref{chemo_4plot}, the Kaplan-Meier (K–M) survival curve, as well as the theoretical and empirical cumulative distribution functions (ECDF), and total time on test (TTT) scaled are displayed. The closely matched empirical and theoretical plots suggest that our model is an excellent fit for the given data. Additionally, the TTT scaled plot shows that the model is suitable for a hazard rate structure that is not monotonic.

\begin{table}[ht]
    \centering
    \caption{MLEs and Goodness-of-Fit Statistics for Chemotherapy Data}
    \label{mle_compare_chemo}
    \scalebox{0.5}{
    \begin{tabular}{llllllllllllll}
    \toprule
& \multicolumn{4}{c}{\textbf{Estimates} (SE)} &  \multicolumn{9}{c}{\textbf{Statistics}} \\ \cline{2-14}

\textbf{Model}&$\alpha$&$\beta$&$\gamma$&&$-2\log \,L$ &$AIC$&$CAIC$&$BIC$&$HQIC$&$W^*$&$A^*$&K-S&p-value\\

\hline

T2GWE&1.1328&0.5416&1.4015 & - & 113.3334 & 119.3334 & 119.9188 & 124.7534 & 121.3539 & 0.0415 & 0.3113 & 0.0756 & 0.9421 \\
&(0.4388)&(0.1170)&(0.5530)\\

\hline

&$\alpha$&$\phi$&$\theta$&&&&&&\\
EGT& 1000.1282 & 0.1452 &  7.1554 & - & 115.9096 & 121.9096 & 122.4949 & 127.3295 & 123.9301 & 0.0608 & 0.4231 & 0.0927 & 0.0926 \\
 &($3.6913\times 10^3$)&($7.4598\times 10^{-2}$)&(3.7373) \\ 

\hline

&$\alpha$&$\theta$&$\gamma$&&&&&&&\\
WGE& 3.9393 & 0.9508 & 0.1484 & - & 115.9251 & 121.9251 & 122.5105 & 127.3451 & 123.9457 & 0.0917 & 0.6079 & 0.1120 &  0.5864 \\
& (8.8328) & (0.2076) & (0.2675)&&&&&&&&\\  

\hline

&$\alpha$&$\beta$&$\theta$&$k$&&&&&\\
LGT & 17.9903 & 0.0196 & 7.0332 & 0.1503 & 116.3564 & 124.3564 & 125.3564 & 131.5831 & 127.0504 & 0.0610 & 0.4270 & 0.0892 & 0.835 \\
& (24.8990) & (0.0298) & (1.9777)&(0.0459)\\ 

\hline

&$\alpha$&$\nu$&\\
T2G&0.4987& 0.8672 & - & - & 127.6381 & 131.6381 & 131.9238 & 135.2515 & 132.9851 & 0.1430 & 0.9790 & 0.1382 & 0.3253 \\
&(0.0979)&(0.0928) \\ 

\hline

&$\alpha$&$\beta$&$\lambda$&$\theta$&\\
EWL & 4.6147 & 0.3981 & 0.4506 & 221.7210 & 140.0032 & 148.0032 & 149.0032 & 155.2298 & 150.6972 & 0.3318 & 2.0776 & 0.1853 & 0.0796 \\
&(1.1065)&(0.6466)&(0.7320)&(250.5359) \\

\bottomrule
\end{tabular}
}
\end{table}

\subsection{Depressive data}

This dataset consists of scores from the ``General Rating of Affective Symptoms for Preschoolers" (GRASP) scale, which is used to assess behavioral and emotional problems in children.
The data set is studied by \cite{LeivaDepressivedata}. Please refer to the ``Declarations" section for detailed instructions on how to access the data.

\begin{figure}[htbp!]
    \centering
    \begin{minipage}{0.45\linewidth}
    \centering
    \includegraphics[width=0.9\linewidth]{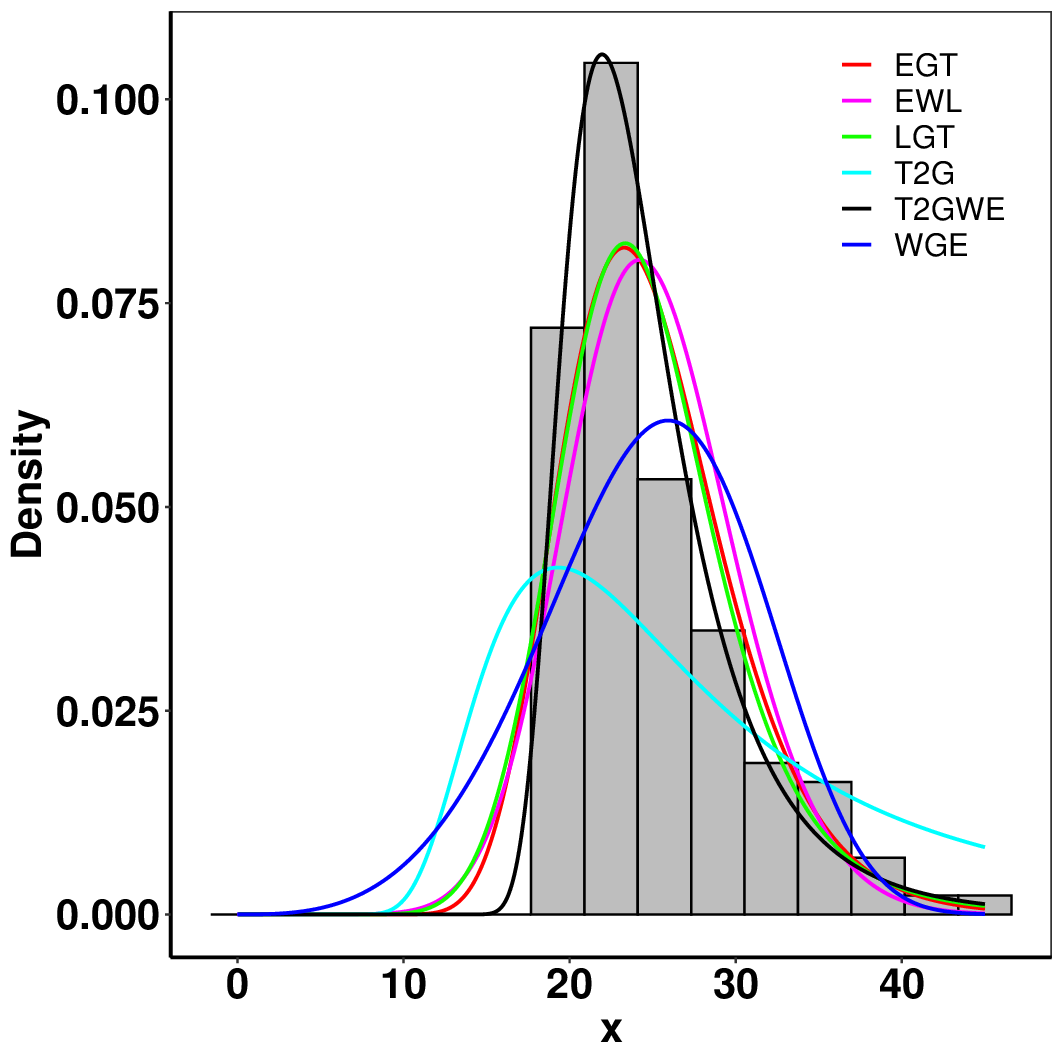}
    \end{minipage}
    \begin{minipage}{0.45\linewidth}
    \centering
    \includegraphics[width=0.9\linewidth]{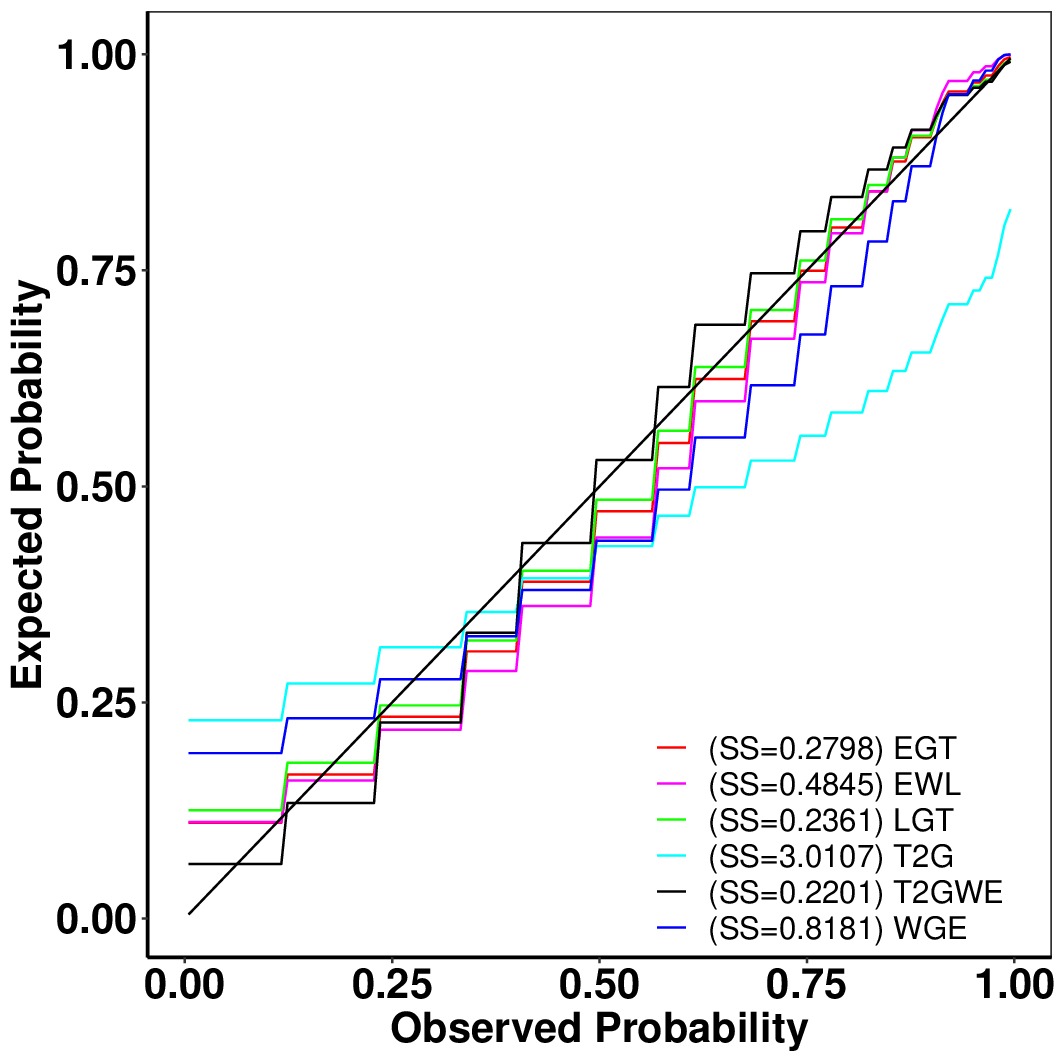}
    \end{minipage}
    \caption{left) Fitted density superposed on the histogram and observed probability for the Depressive data. right) Expected probability plots for the Depressive data.}
    \label{hist_ep_depressive}
\end{figure}
\begin{figure}[htbp!]
    \centering
    \begin{minipage}{0.45\linewidth}
    \centering
    \includegraphics[width=0.9\linewidth]{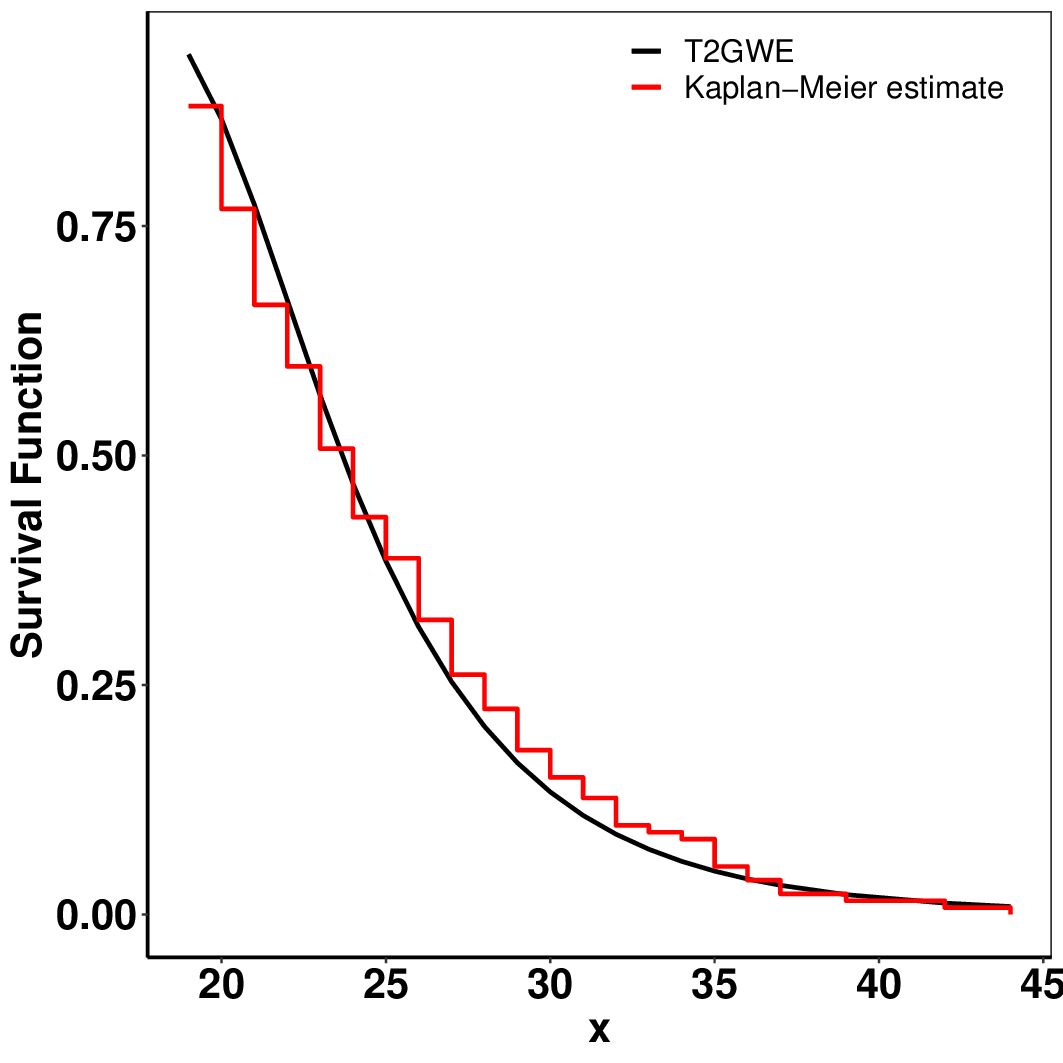}
    \end{minipage}
    \begin{minipage}{0.45\linewidth}
    \centering
    \includegraphics[width=0.9\linewidth]{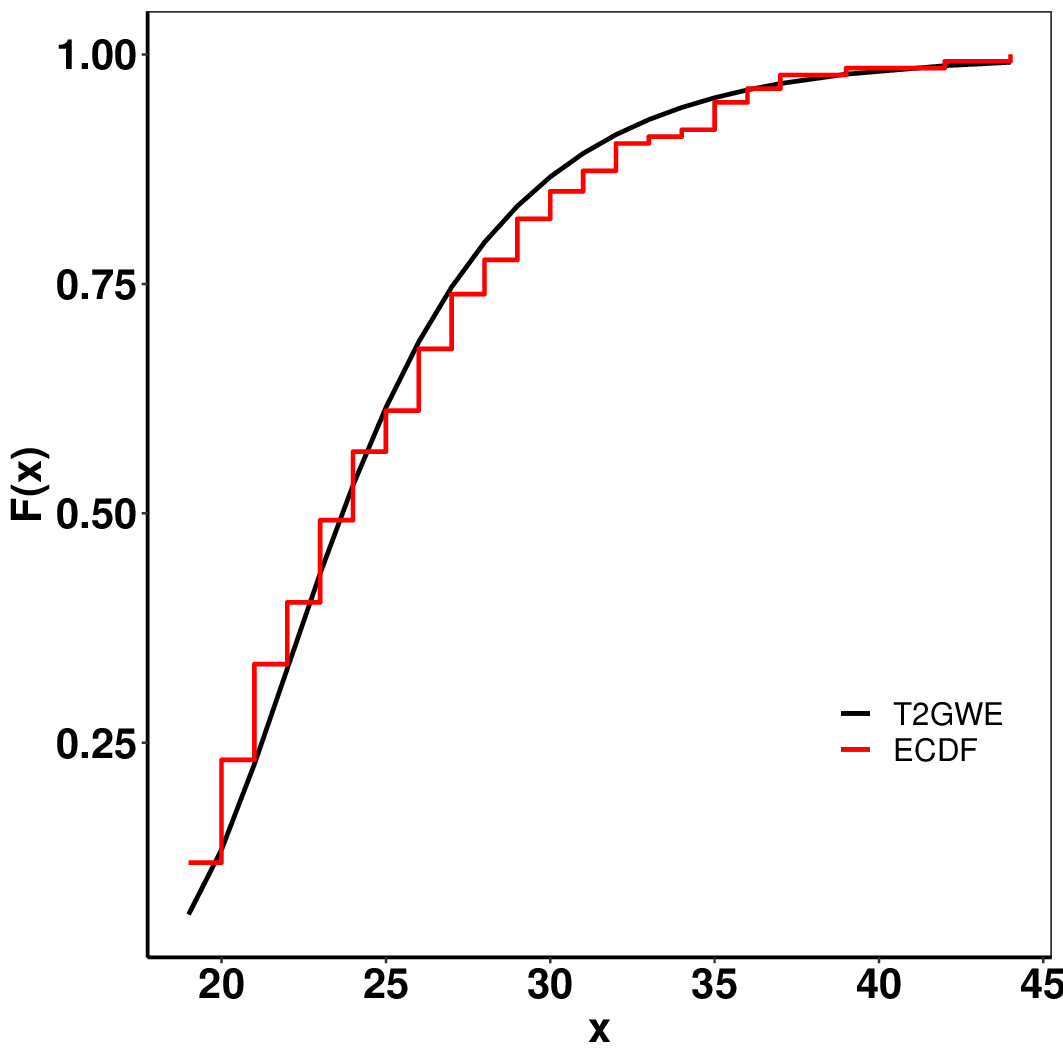}
    \end{minipage}
    \qquad
    
    \begin{minipage}{0.45\linewidth}
    \centering
    \includegraphics[width=0.9\linewidth]{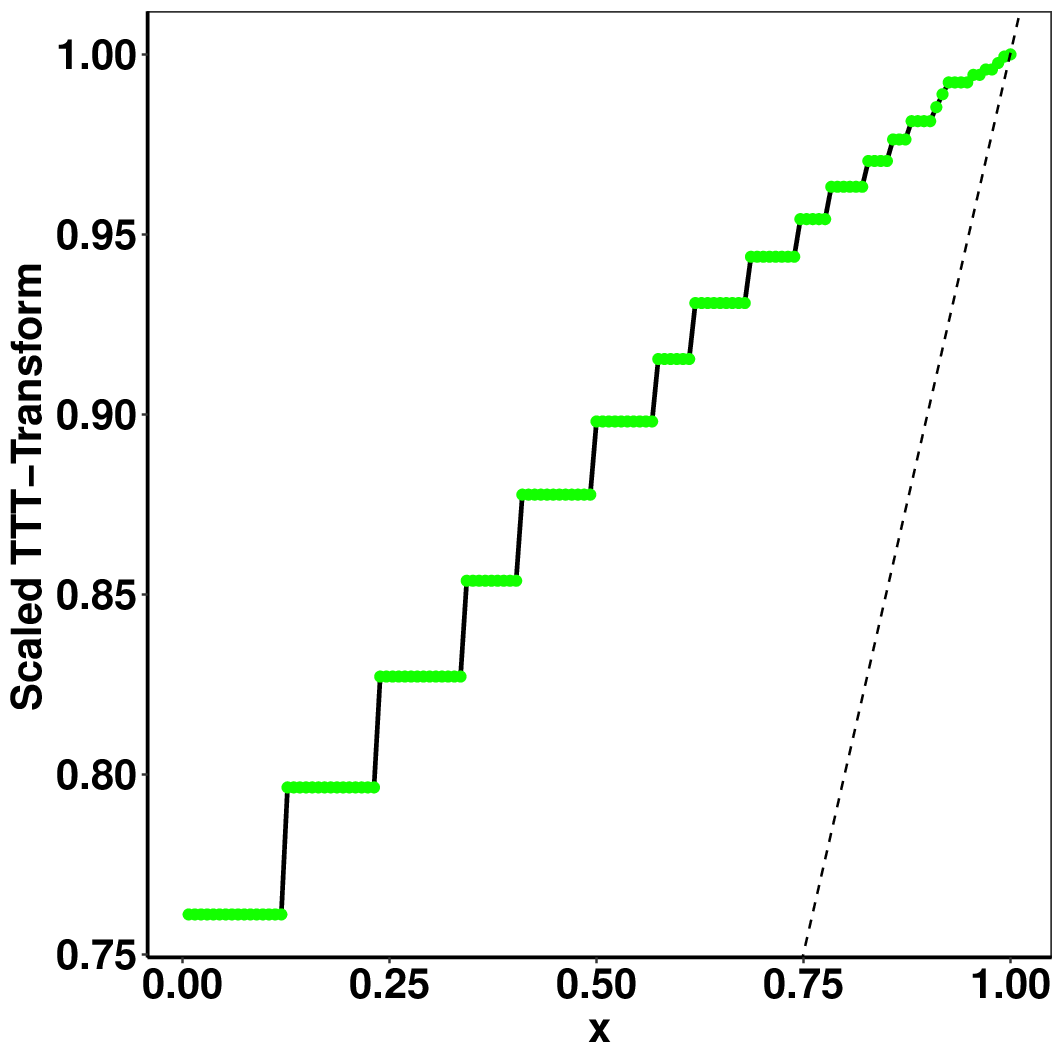}
    \end{minipage}
    \begin{minipage}{0.45\linewidth}
    \centering
    \includegraphics[width=0.9\linewidth]{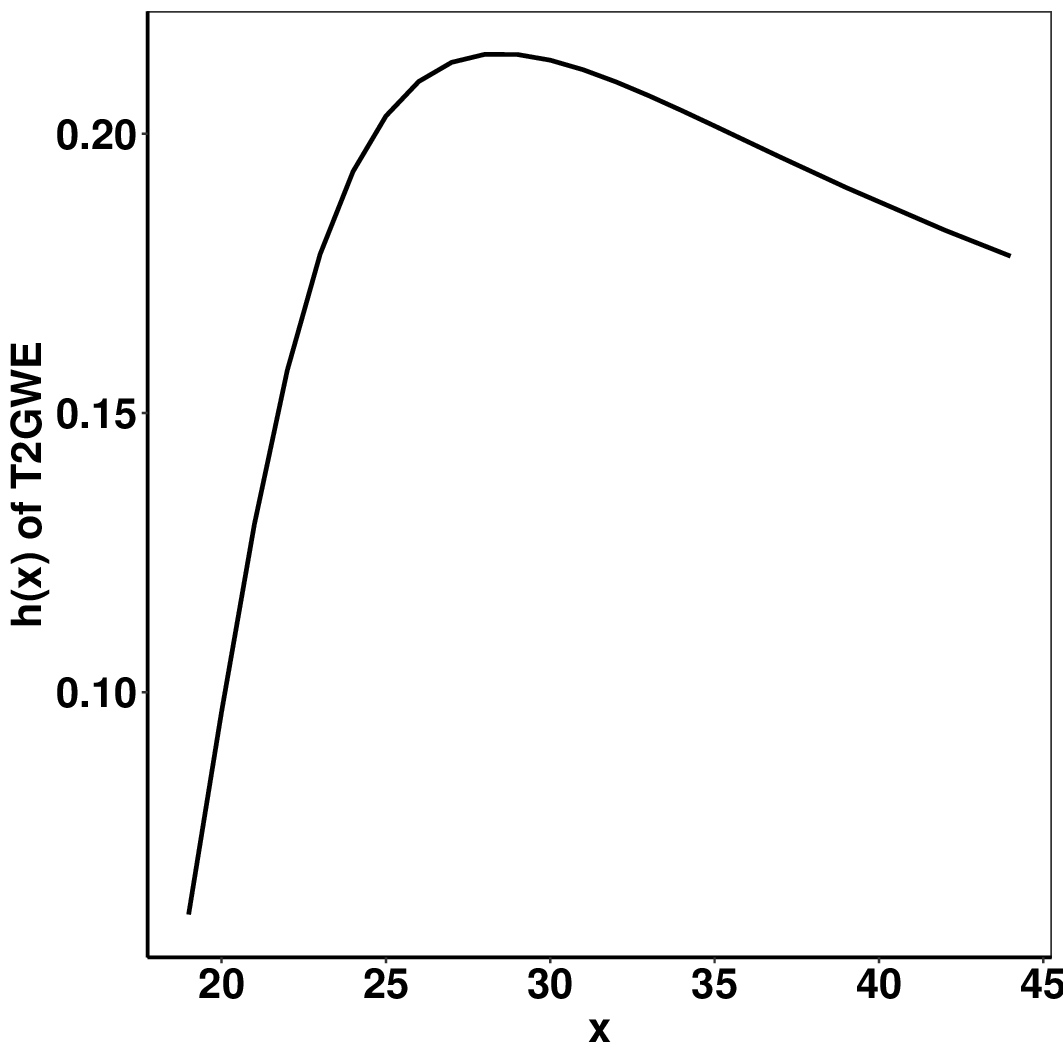}
    \end{minipage}
    \caption{Fitted K-M survival curve, theoretical and ECDF, the TTT statistics, and the hrf for the Depressive data.}
    \label{depressive_4plot}
\end{figure}

From Table \ref{mle_compare_depressive} and Fig.~\ref{hist_ep_depressive}, we can conclude that the Type-2 Gumbel Weibull-exponential distribution has better performance than other distributions as it possesses the smallest values for all goodness-of-fit statistics and the highest p-value in the K-S test. As shown in Fig.~\ref{depressive_4plot}, the close resemblance between the fitted empirical and theoretical plots indicates a strong fit of our model to the provided data.
Furthermore, the TTT scaled plot provides clear evidence that the model is appropriate for a non-monotonic hazard rate structure.

\begin{table}[ht]
    \centering
    \caption{MLEs and Goodness-of-Fit Statistics for Depressive Data}
    \label{mle_compare_depressive}
    \scalebox{0.5}{
    \begin{tabular}{llllllllllllll}
    \toprule
& \multicolumn{4}{c}{\textbf{Estimates} (SE)} &  \multicolumn{9}{c}{\textbf{Statistics}} \\ \cline{2-14} \textbf{Model} &$\alpha$&$\beta$&$\gamma$&&$-2\log \,L$ &$AIC$&$CAIC$&$BIC$&$HQIC$&$W^*$&$A^*$&K-S&p-value\\

\hline

T2GWE & 0.1127 & 5.0122 & 0.0223 & - & 785.7463 & 791.7463 & 791.9309 & 800.4398 & 795.2791 & 0.1814 & 1.3478 & 0.1092 & 0.0821 \\
&(0.0818)&(0.3542)&(0.0024) \\

\hline

&$\alpha$&$\phi$&$\theta$&&&&&&\\
EGT& 9.7621 & 2.0192 & 1691.0601 & - & 803.2078 & 809.2078 & 809.3924 & 817.9013 & 812.7406 & 0.3185 & 2.1927 & 0.1108 & 0.0746 \\
&(3.7163)&(0.2403)&(1088.5591) \\ 

\hline

&$\alpha$&$\theta$&$\gamma$&&&&&&&\\
WGE& 44.2733 & 3.7778 & 0.0115 & - & 852.253 & 858.253 & 858.4376 & 866.9465 & 861.7858 & 0.8203 & 5.0577 & 0.1913 & 0.0001 \\
& (163.1340) & (0.5147) &(0.0008) &&&&&&&&\\  

\hline

&$\alpha$&$\beta$&$\theta$&$k$&&&&&\\
LGT & 2.7374 & 0.0006 & 161.2461 & 1.0193 & 807.7372 & 815.7372 & 816.0473 & 827.3286 & 820.4476 & 0.3392 & 2.3111 & 0.1252 & 0.0301 \\
& (0.9395) & (0.0028) & (46.6043)&(0.1060)\\ 

\hline

&$\alpha$&$\nu$&\\
T2G& 1684.3935 & 2.3918 & - & - & 919.7708 & 923.7708 & 923.8624 & 929.5665 & 926.126 & 0.2340 & 1.6765 & 0.2480 & $1.392e-07$ \\
&(404.6525)&(0.0807) \\ 

\hline

&$\alpha$&$\beta$&$\lambda$&$\theta$&\\
EWL & 1.2679 & 0.0770 & 0.6636 & 61.1849 & 816.1935 & 824.1935 & 824.5036 & 835.7849 & 828.9039 & 0.4525 & 2.9886 & 0.1310 & 0.0202 \\
&(0.3130)&(0.0271)&(0.2349)&(31.6238)\\

\bottomrule
\end{tabular}
}
\end{table}

\subsection{Covid-Mexico data}
This dataset recorded the mortality rates of the patients infected by the COVID-19 pandemic in Mexico and was studied by \cite{zhou2023implementation}. This dataset consists of 106 observations from March 31, 2020, to July 20, 2020. 
To access the data, please follow the instructions in the ``Declarations" section.

\begin{figure}[htbp!]
    \centering
    \begin{minipage}{0.45\linewidth}
    \centering
    \includegraphics[width=0.9\linewidth]{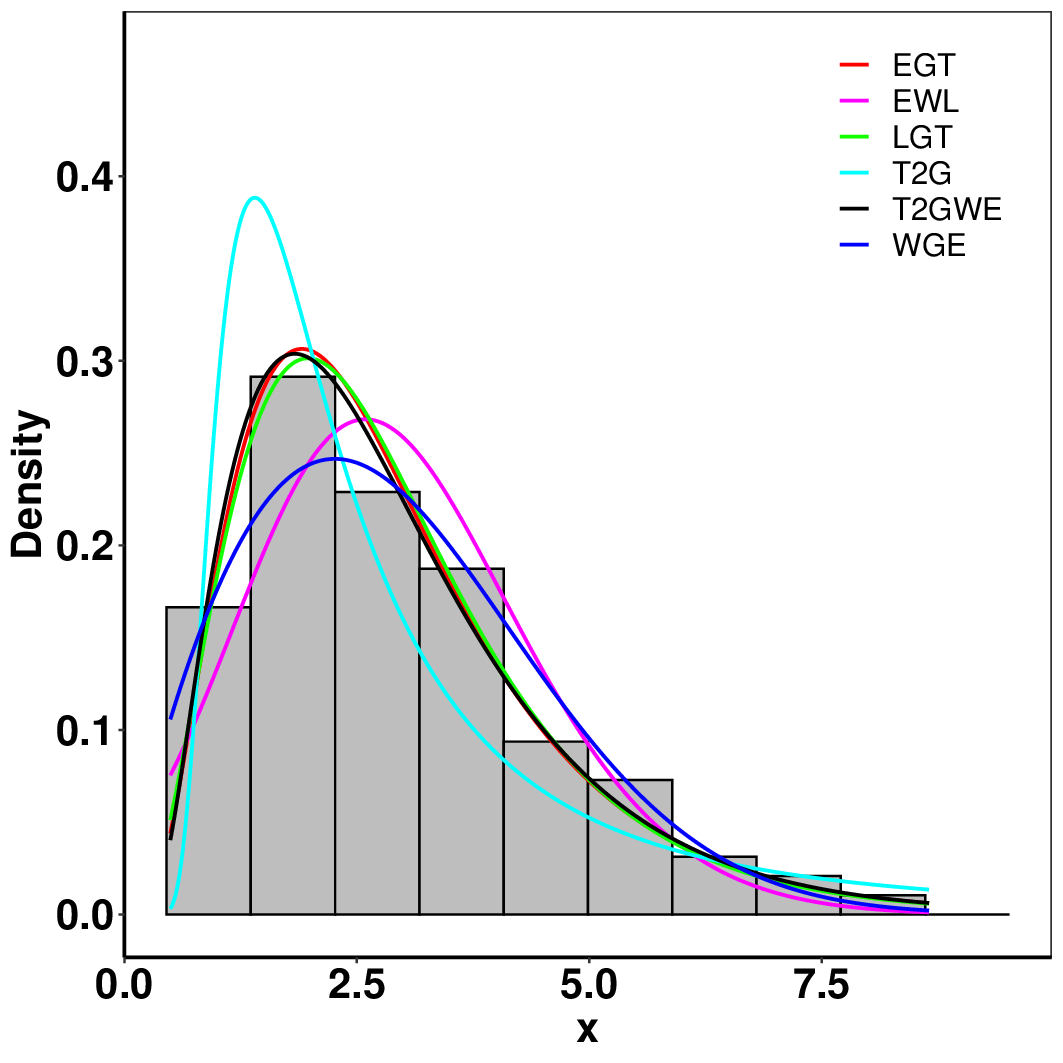}
    \end{minipage}
    \begin{minipage}{0.45\linewidth}
    \centering
    \includegraphics[width=0.9\linewidth]{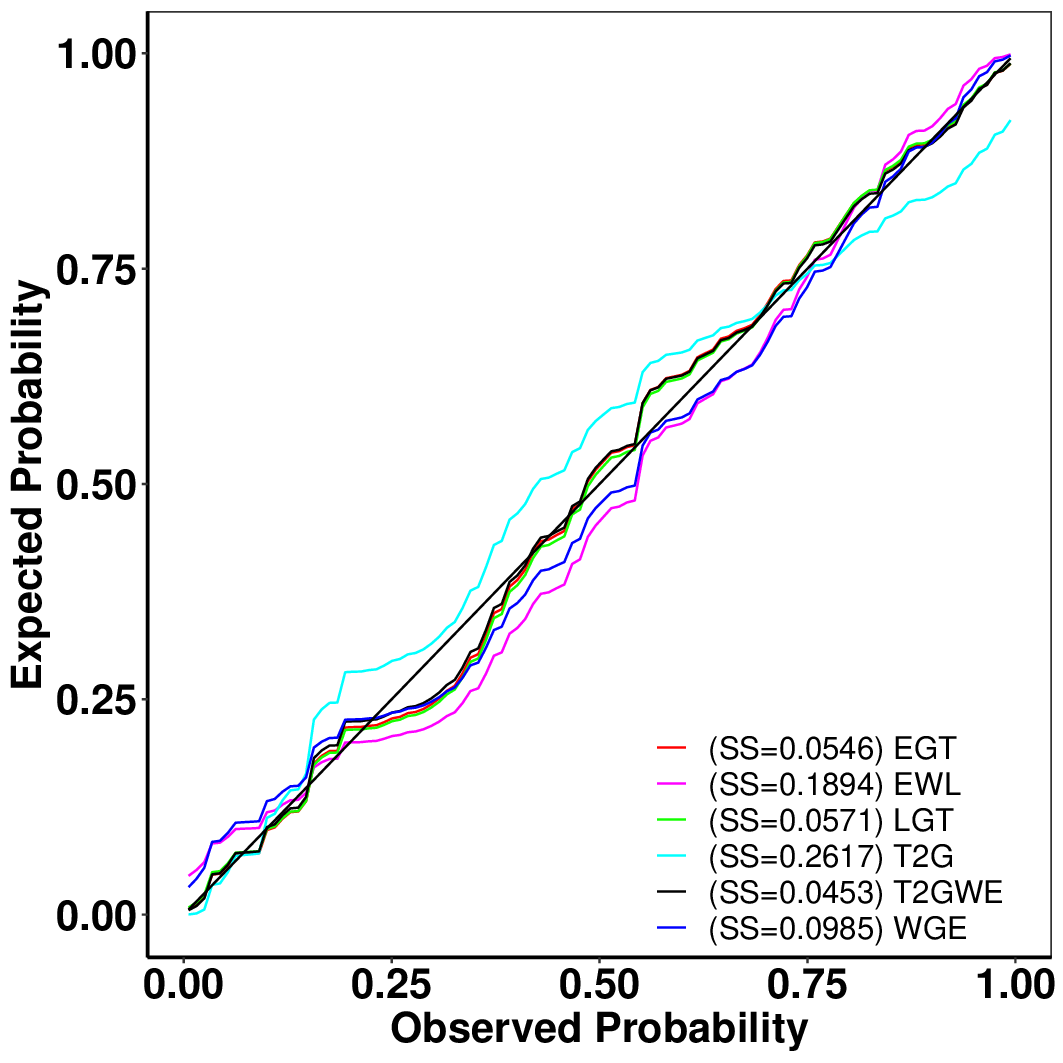}
    \end{minipage}
    \caption{left) Fitted density superposed on the histogram and observed probability for the Covid-Mexico data. right) Expected probability plots for the Covid-Mexico data.}
    \label{hist_ep_mexico}
\end{figure}

\begin{figure}[htbp!]
    \centering
    \begin{minipage}{0.45\linewidth}
    \centering
    \includegraphics[width=0.9\linewidth]{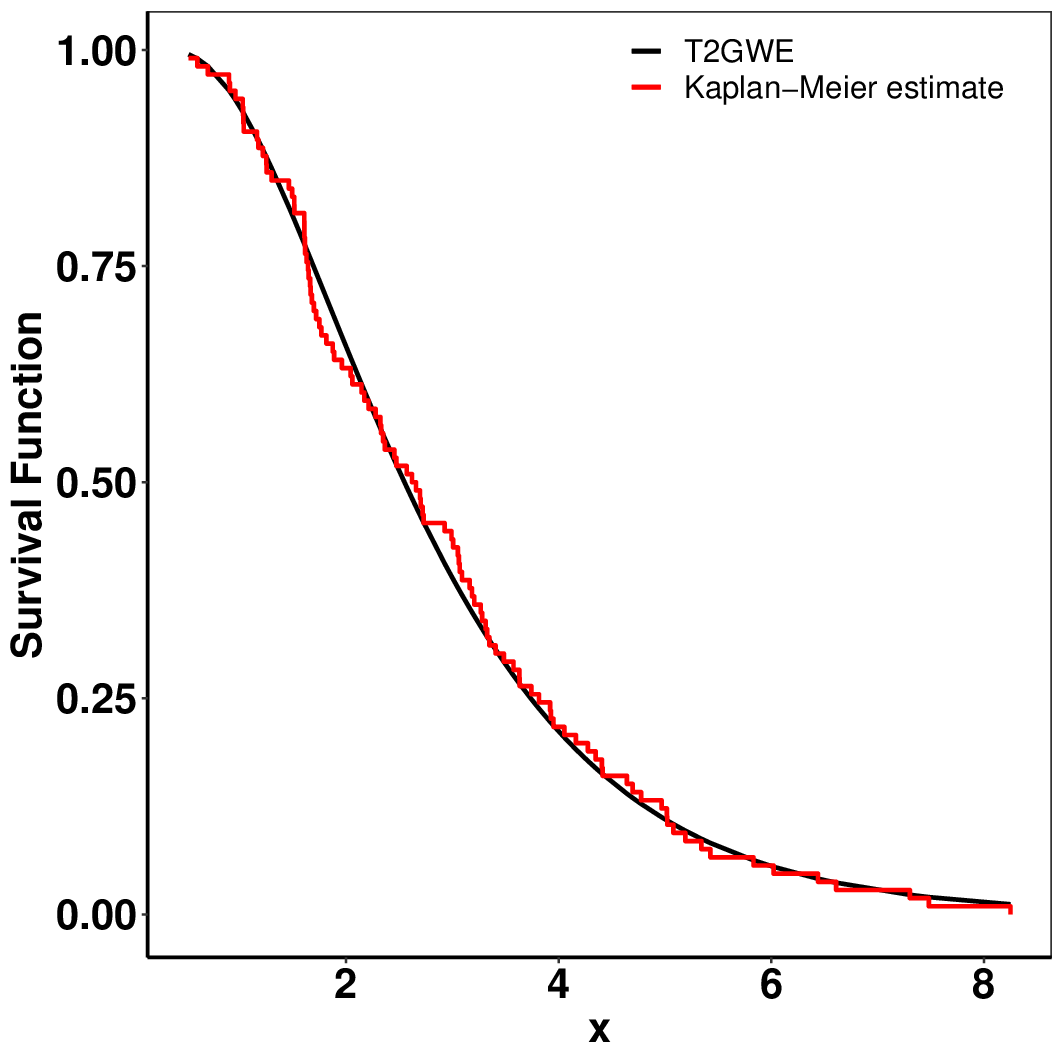}
    \end{minipage}
    \begin{minipage}{0.45\linewidth}
    \centering
    \includegraphics[width=0.9\linewidth]{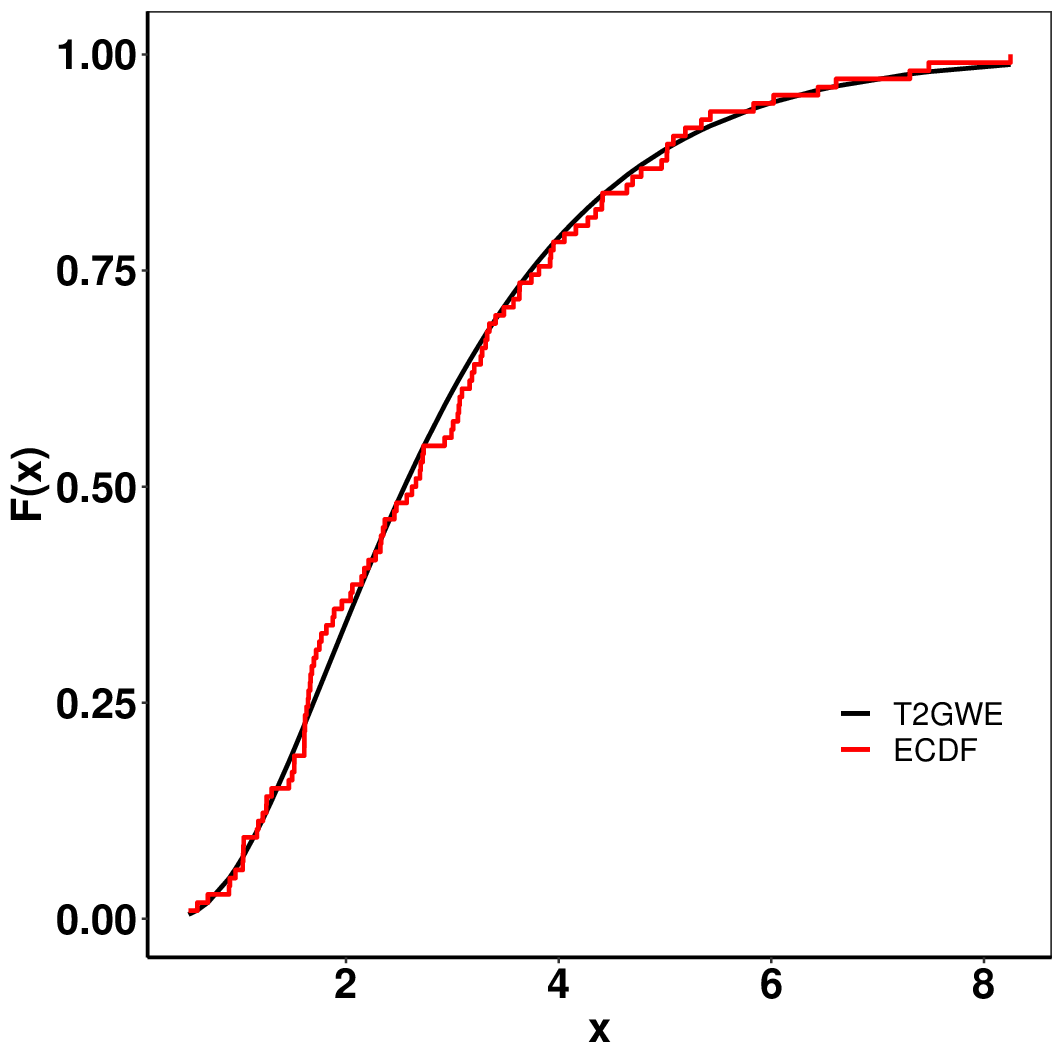}
    \end{minipage}
    \qquad
    
    \begin{minipage}{0.45\linewidth}
    \centering
    \includegraphics[width=0.9\linewidth]{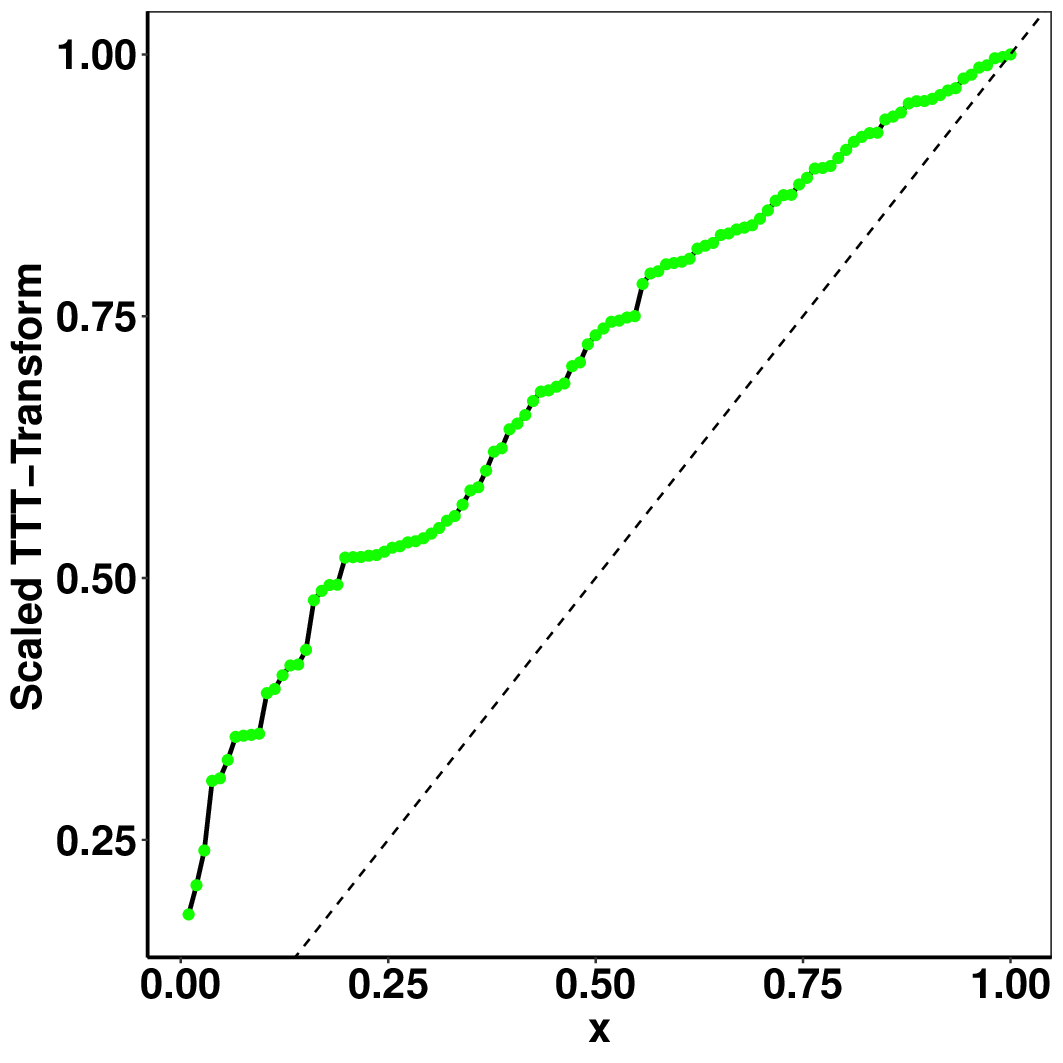}
    \end{minipage}
    \begin{minipage}{0.45\linewidth}
    \centering
    \includegraphics[width=0.9\linewidth]{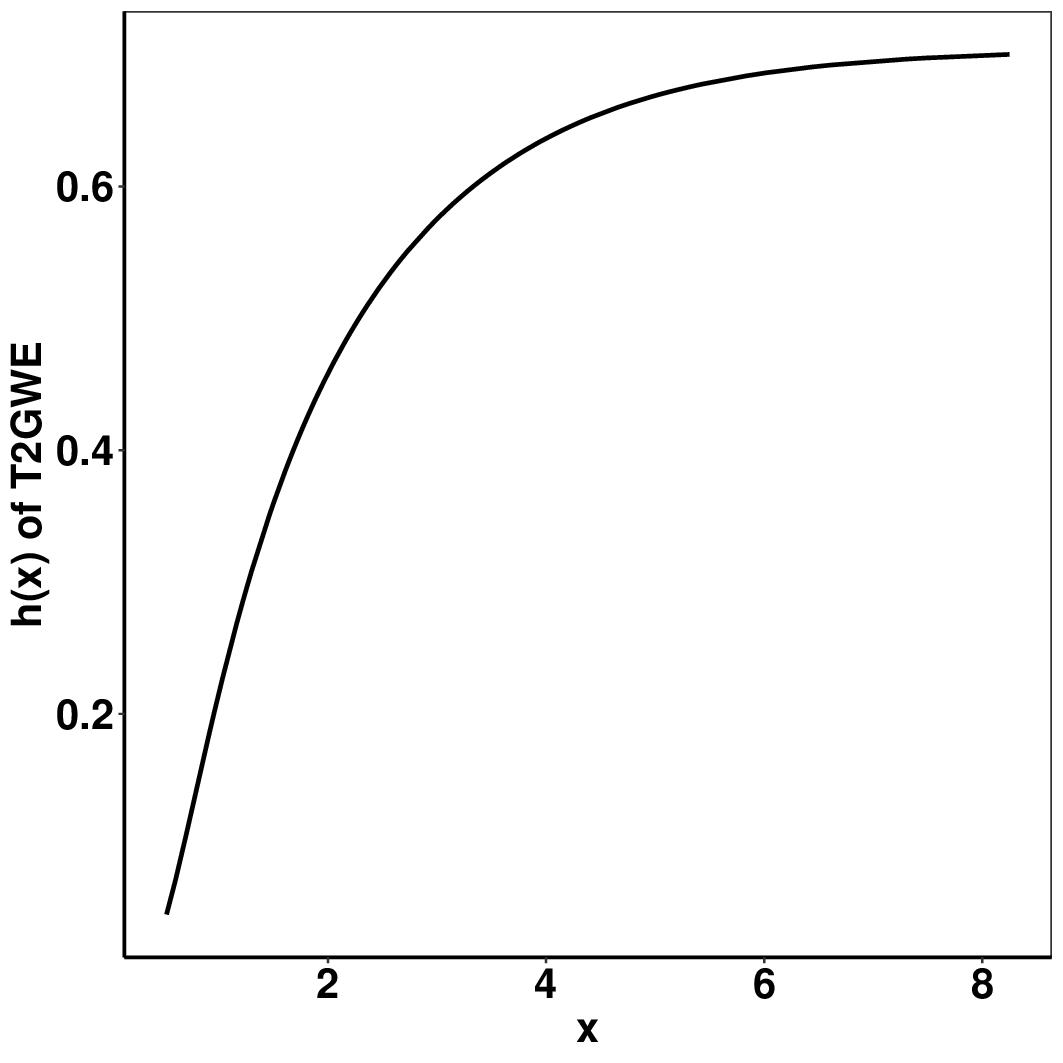}
    \end{minipage}
    \caption{Fitted K-M survival curve, theoretical and ECDF functions, the TTT statistics, and the hrf for the Covid-Mexico data.}
    \label{mexico_4plot}
\end{figure}

Table \ref{mle_compare_covid_mexico} and Fig.~\ref{hist_ep_mexico} provide
the estimates of the parameters and the goodness-of-fit.  We can conclude that the Type-2 Gumbel Weibull-exponential distribution has the best performance among other distributions given its smallest values in all goodness-of-fit statistics and the highest p-value in the K-S test.
Displayed in Fig. \ref{mexico_4plot} are the K-M survival curve, as well as the theoretical and ECDF, and TTT scaled. The convergence of the fitted empirical and theoretical plots suggests that our model accurately represents the given data.
In addition, the TTT scaled plot clearly indicates that the model is well-suited for a hazard rate structure that is not strictly monotonic.

\begin{table}[ht]
    \centering
    \caption{MLEs and Goodness-of-Fit Statistics for Covid-Mexico Data}
    \label{mle_compare_covid_mexico}
    \scalebox{0.5}{
    \begin{tabular}{llllllllllllll}
    \toprule
& \multicolumn{4}{c}{\textbf{Estimates} (SE)} &  \multicolumn{9}{c}{\textbf{Statistics}} \\\cline{2-14}
\textbf{Model}&$\alpha$&$\beta$&$\gamma$&&$-2\log \,L$ &$AIC$&$CAIC$&$BIC$&$HQIC$&$W^*$&$A^*$&K-S&p-value\\

\hline

T2GWE & 3.9220 & 0.7209 & 0.9767 & - & 375.7089 & 381.7089 & 381.9442 & 389.6993 & 384.9475 & 0.0454 & 0.2525 & 0.0577 & 0.8721 \\
&(1.4023)&(0.2117)&(0.3942)\\

\hline

&$\alpha$&$\phi$&$\theta$&&&&&&\\
EGT& 56.8310 & 0.4537 & 6.7577 & - & 376.3527 & 382.3527 & 382.5879 & 390.343 & 385.5912 & 0.0521 & 0.2877 & 0.0646 & 0.7686 \\
&(78.7004)&(0.1409)&(1.3261) \\ 

\hline

&$\alpha$&$\theta$&$\gamma$&&&&&&&\\
WGE& 49.8868 & 1.8037 & 0.0326 & - & 382.8119 & 388.8119 & 389.0472 & 396.8022 & 392.0504 & 0.1128 & 0.7313 & 0.0664 & 0.7387 \\
& (72.3253) & (0.1413) &(0.0222) & \\  

\hline

&$\alpha$&$\beta$&$\theta$&$k$&&&&&\\
LGT & 11.9707 & 0.0343 & 8.5300 & 0.3374 & 376.8091 & 384.8091 & 385.2052 & 395.4629 & 389.1271 & 0.0557 & 0.3106 & 0.0688 & 0.6975 \\
& (11.5930) & (0.0615) & (1.6317)&(0.0982)\\ 

\hline

&$\alpha$&$\nu$&\\
T2G& 2.8185 & 1.6830 & - & - & 400.2692 & 404.2692 & 404.3857 & 409.5961 & 406.4282 & 0.2734 & 1.7490 & 0.0928 & 0.3213 \\
&(0.2928)&(0.1163) \\ 

\hline

&$\alpha$&$\beta$&$\lambda$&$\theta$&\\
EWL & 3.6179 & 0.3384 & 0.4196 & 150.8827 & 390.9746 & 398.9747 & 399.3707 & 409.6284 & 403.2927 & 0.1807 & 1.1892 & 0.0957 & 0.2857 \\
&(0.6784)&(2.5979)&(3.2208)&(109.2607)\\

\bottomrule
\end{tabular}
}
\end{table}

\section{Conclusion}

This paper introduces a novel methodology for generating continuous statistical distributions by utilizing the exponentiated odds ratio, which is based on the concepts of survival analysis. The approach described herein constitutes a substantial advancement in the field of statistical modeling, effectively augmenting the flexibility and accuracy of distribution models to effectively address the requirements posed by intricate contemporary data architectures. The major point of this progress lies in the formulation of the ``Type-2 Gumbel Weibull-G Family of Distributions," which has undergone a comprehensive mathematical analysis. The scope of this investigation covered various statistical properties, such as expansions of density functions, moments, hazard rate and quantile functions, R\'enyi entropy, order statistics, and an examination of stochastic ordering.

In order to assess the robustness and reliability of the new generator, we employed a set of five advanced parameter estimation techniques: Maximum Likelihood, Least Square, Weighted Least Square, Maximum Product Spacing, Cram\'er-von Mises, and Anderson and Darling. The efficacy and utility of the Type-2 Gumbel Weibull-G distributions were further validated through a comprehensive analysis of three datasets obtained from real-world scenarios. These practical implementations demonstrated the superior statistical accuracy of our proposed distributions over existing models, thereby emphasizing their relevance and applicability in both theoretical and practical statistical domains.

Our future research endeavors will entail a thorough investigation of several sub-families within the new generator. This exploration will focus on their distinct properties and potential applications in diverse scenarios, setting them in comparison with other established distribution models. Furthermore, we are now developing an R package with the objective of optimizing the parameter estimate process through the utilization of diverse methodologies. This will result in improved efficiency of data-fitting procedures, thus increasing the accessibility and practicality of our study for wider applications.

\section*{Abbreviations}{
The following abbreviations are used in this paper:\\

\noindent 
\begin{tabular}{@{}ll}
T2GWG & The Type-2 Gumbel Weibull-G \\
cdf & cumulative distribution function\\
pdf  & probability density function \\
hrf & hazard rate function \\
Exp-G & exponentiated-G \\
EGT & Exponentiated Gumbel Type-2 \\
WGE & Weibull Generalized Exponential \\
LGT & Lomax Gumbel Type-2 \\
T2G & Type-2 Gumbel \\
EWL & Exponentiated Weibull-Logistic distributions \\
MLE & maximum likelihood estimates \\
MPS & maximum product spacing estimates \\
LS & least square estimates \\
WLS & weighted least square estimates \\
CVM & Cram\'er-von Mises estimates \\
AD & Anderson and Darling estimates \\
T2GWE & Type-2 Gumbel Weibull-Exponential \\
T2GWU & Type-2 Gumbel Weibull-Uniform \\
T2GWP & Type-2 Gumbel Weibull-Pareto \\
AIC & Akaike Information Criterion \\
CAIC & Consistent Akaike Information Criterion \\
BIC & Bayesian Information Criterion \\
HQIC & Hannan-Quinn Criterion \\
$W^{*}$ & Cram\'er-von Mises statistic \\
$A^{*}$  & Anderson-Darling statistic\\
K-S & Kolmogorov-Smirnov statistic \\
ECDF & empirical cumulative distribution function \\
TTT & total time on test\\
K-M & Kaplan-Meier
\end{tabular}
}

\backmatter

\bmhead{Supplementary information}
Not applicable


\bmhead{Acknowledgments}

The authors wish to thank all members of the CSDA lab at the University of West Florida for their helpful comments on the manuscript.
\section*{Declarations}

\begin{itemize}
\item Funding: This research received no specific grant from any funding agency in the public, commercial, or
not-for-profit sectors.
\item Conflict of interest/Competing interests: On behalf of all authors, the corresponding author states that there is no conflict of interest.
\item Ethics approval: Not applicable
\item Consent to participate: Not applicable
\item Consent for publication: Not applicable
\item Availability of data and materials: All data utilized in this study is openly available on GitHub. The dataset can be accessed through the provided link: \url{https://github.com/shusenpu/Generator_Data}. The data can also be accessed using its Digital Object Identifier (DOI): 10.5281/zenodo.10215787. To retrieve the data, visit the following link: \url{https://doi.org/10.5281/zenodo.10215787}.
\item Code availability: available upon request.
\item Authors' contributions: All authors have contributed significantly to writing and editing the paper.
\end{itemize}

\noindent






\begin{appendices}

\section{The first derivatives of H}\label{sec:MPS_pd}
The first partial derivatives of $H$ in Section 4 with respect to $\alpha, \beta, \delta, \psi_k$ are given by 
\begin{align}
    \frac{\partial H}{\partial\alpha}=& \frac{1}{n+1}\left\{-\left[ \frac{H(x_1,\bpsi)}{\overline{H}(x_1,\bpsi)} \right]^{-\beta} +\frac{\left[\frac{H(x_n,\bpsi)}{\overline{H}(x_n,\bpsi)} \right]^{-\beta}\exp\left\{-\alpha \left[ \frac{H(x_i,\bpsi)}{\overline{H}(x_i,\bpsi)} \right]^{-\beta} \right\}}{1-\exp\left\{-\alpha \left[ \frac{H(x_i,\bpsi)}{\overline{H}(x_i,\bpsi)} \right]^{-\beta} \right\}} \right. \nonumber \\
    &+\sum_{i=2}^n \left[\frac{\left[ \frac{H(x_i,\bpsi)}{\overline{H}(x_i,\bpsi)} \right]^{-\beta} \exp\left\{-\alpha\left[ \frac{H(x_i,\bpsi)}{\overline{H}(x_i,\bpsi)} \right]^{-\beta} \right\}}{\exp\left\{-\alpha\left[ \frac{H(x_i,\bpsi)}{\overline{H}(x_i,\bpsi)} \right]^{-\beta} \right\}-\exp\left\{-\alpha\left[ \frac{H(x_{i-1},\bpsi)}{\overline{H}(x_{i-1},\bpsi)} \right]^{-\beta} \right\}} \right. \nonumber \\
    & \left. \left.-\frac{\left[ \frac{H(x_{i-1},\bpsi)}{\overline{H}(x_{i-1},\bpsi)} \right]^{-\beta}\exp \left\{-\alpha \left[ \frac{H(x_{i-1},\bpsi)}{\overline{H}(x_{i-1}, \bpsi)} \right]^{-\beta} \right\}}{\exp\left\{-\alpha\left[ \frac{H(x_i,\bpsi)}{\overline{H}(x_i,\bpsi)} \right]^{-\beta} \right\}-\exp\left\{-\alpha\left[ \frac{H(x_{i-1},\bpsi)}{\overline{H}(x_{i-1},\bpsi)} \right]^{-\beta} \right\}} \right] \right\}
\end{align}
\begin{align}
    \frac{\partial H}{\partial\beta}=& \frac{1}{n+1}\left\{\alpha\log \left(\frac{H(x_1,\bpsi)}{\overline{H}(x_1,\bpsi)}\right) \left[\frac{H(x_i,\bpsi)}{\overline{H}(x_i,\bpsi)} \right]^{-\beta} \right. \nonumber\\
    &+ \frac{\log\left( \frac{H(x_n,\bpsi)}{\overline{H}(x_n,\bpsi)}\right) \left[ \frac{H(x_n,\bpsi)}{\overline{H}(x_n,\bpsi)} \right]^{-\beta} \exp\left\{- \alpha \left[ \frac{H(x_i,\bpsi)}{\overline{H}(x_i,\bpsi)} \right]^{-\beta} \right\}}{1-\exp\left\{-\alpha \left[\frac{H(x_i,\bpsi)}{\overline{H}(x_i,\bpsi)} \right]^{-\beta} \right\}} \nonumber \\
    &+\sum_{i=2}^n  \left[ \frac{\log\left( \frac{H(x_i,\bpsi)}{\overline{H}(x_i,\bpsi)}\right) \left[ \frac{H(x_i,\bpsi)}{\overline{H}(x_i,\bpsi)} \right]^{-\beta} \exp \left\{-\alpha\left[ \frac{H(x_i, \bpsi)}{\overline{H}(x_i,\bpsi)} \right]^{-\beta} \right\}}{\exp\left\{-\alpha \left[ \frac{H(x_i,\bpsi)}{\overline{H}(x_i,\bpsi)} \right]^{-\beta} \right\}-\exp\left\{-\alpha\left[ \frac{H(x_{i-1},\bpsi)}{\overline{H}(x_{i-1},\bpsi)} \right]^{-\beta} \right\}} \right. \nonumber \\
    & \left. \left. \left.- \frac{\log\left( \frac{H(x_{i-1},\bpsi)}{\overline{H}(x_{i-1},\bpsi)}\right) \left[ \frac{H(x_{i-1},\bpsi)}{\overline{H}(x_{i-1},\bpsi)} \right]^{-\beta}\exp\left\{-\alpha\left[ \frac{H(x_{i-1}, \bpsi)}{\overline{H}(x_{i-1}, \bpsi)} \right]^{-\beta} \right\}}{\exp\left\{-\alpha \left[ \frac{H(x_i,\bpsi)}{\overline{H}(x_i,\bpsi)} \right]^{-\beta} \right\}-\exp\left\{-\alpha\left[ \frac{H(x_{i-1},\bpsi)}{\overline{H}(x_{i-1},\bpsi)} \right]^{-\beta} \right\}} \right] \right) \right\}
\end{align}

\begin{align}
    \frac{\partial H}{\partial\psi_k}=& \frac{1}{n+1} \left\{\alpha\beta \frac{H(x_1,\bpsi)^{-\beta-1}}{\overline{H}(x_1,\bpsi)^{-\beta+1}} \frac{\partial H(x_1,\bpsi)}{\partial\psi_k} + \frac{\frac{H(x_n,\bpsi)^{-\beta-1}}{\overline{H}(x_n,\bpsi)^{-\beta+1}} \frac{\partial H(x_n,\bpsi)}{\partial\psi_k} \exp\left\{-\alpha \left[ \frac{H(x_n,\bpsi)}{\overline{H}(x_n,\bpsi)} \right]^{-\beta} \right\}}{1-\exp\left\{-\alpha \left[ \frac{H(x_n,\bpsi)}{\overline{H}(x_n,\bpsi)} \right]^{-\beta} \right\}} \right.  \nonumber \\
    &+\sum_{i=2}^n \left[ \frac{ \frac{H(x_i,\bpsi)^{-\beta-1}}{\overline{H}(x_i,\bpsi)^{-\beta+1}}\frac{\partial H(x_i,\bpsi)}{\partial\psi_k} \exp\left\{-\alpha \left[ \frac{H(x_i,\bpsi)}{\overline{H}(x_i,\bpsi)}\right]^{-\beta} \right\}}{\exp\left\{-\alpha\left[ \frac{H(x_i,\bpsi)}{\overline{H}(x_i,\bpsi)} \right]^{-\beta} \right\}-\exp\left\{-\alpha \left[\frac{H(x_{i-1},\bpsi)}{\overline{H}(x_{i-1},\bpsi)} \right]^{-\beta}\right\}} \right.\nonumber \\
    & \left. \left. \left.- \frac{\frac{H(x_{i-1},\bpsi)^{-\beta-1}}{\overline{H}(x_{i-1},\bpsi)^{-\beta+1}} \frac{\partial H(x_{i-1},\bpsi)}{\partial\psi_k}\exp\left\{-\alpha\left[\frac{H(x_{i-1},\bpsi)}{\overline{H}(x_{i-1},\bpsi)} \right]^{-\beta} \right\}}{\exp\left\{-\alpha\left[ \frac{H(x_i,\bpsi)}{\overline{H}(x_i,\bpsi)} \right]^{-\beta} \right\}-\exp\left\{-\alpha \left[ \frac{H(x_{i-1},\bpsi)}{\overline{H}(x_{i-1},\bpsi)} \right]^{-\beta} \right\}} \right] \right) \right\}
\end{align}

\section{Distributions Used in the Application Section}\label{app:distributions}
\begin{itemize}
    \item Exponentiated Gumbel Type-2 Distribution: 
$$
F_{_{EGT}}(x)=1-\left(1-e^{-\theta x^{-\phi}}\right)^{\alpha},
$$    
$$
f_{_{EGT}}(x)=\alpha \phi \theta x^{-\phi-1} e^{-\theta x^{-\phi}}\left(1-e^{-\theta x^{-\phi}}\right)^{\alpha-1},
$$
$$\alpha, \phi, \theta>0$$
\item Weibull Generalized Exponential Distribution: 
$$
F_{_{WGE}}(x)= 1-e^{-a\left[e^{\lambda x}-1\right]^b},
$$
$$
f_{_{WGE}}(x)=a b \lambda e^{\lambda x}\left[e^{\lambda x}-1\right]^{b-1} e^{-a\left[e^{\lambda x}-1\right]^b},
$$
$$a, b, \lambda>0$$
   \item Lomax Gumbel Type-2 Distribution: 
$$
F_{_{LGT}}(x)=1-\beta^{\alpha}\left( \beta-\log \left[1-e^{-\theta x^{-k}}\right]\right)^{-\alpha}
$$   
$$
f_{_{LGT}}(x)=\frac{\alpha \beta^\alpha \theta k x^{-k-1} e^{-\theta x^{-k}}}{\left[1-e^{-\theta x^{-k}}\right]}\left[\beta-\log \left[1-e^{-\theta x^{-k}}\right]\right]^{-(\alpha+1)}
$$
$$\alpha, \beta, \theta, k>0$$
\item Type-2 Gumbel Distribution:
$$
F_{_{T2G}}(x)=e^{-\theta x^{-\phi}}
$$
$$
f_{_{T2G}}(x)=\phi \theta x^{-\phi-1} e^{-\theta x^{-\phi}}
$$
$$\phi, \theta>0$$
\item Exponentiated Weibull-Logistic Distribution: 
$$
F_{_{EWL}}(x)=\left(1-e^{-\alpha e^{\lambda \beta x}}\right)^{\theta}
$$
$$
f_{_{EWL}}(x)=\theta\left(1-e^{-\alpha e^{\lambda \beta x}}\right)^{\theta-1}\left(\lambda \alpha \beta e^{\lambda \beta x-\alpha e^{\lambda \beta x}}\right),
$$
$$\lambda, \alpha, \beta, \theta>0.$$

\end{itemize}




\end{appendices}


\bibliography{sn-bibliography}

\end{document}